\documentclass[10pt]{amsart}
\usepackage{amsmath}
\usepackage{amsthm}
\usepackage{amscd,amssymb}
\usepackage[all]{xy}
\usepackage{color}
\usepackage{appendix} 
\usepackage{verbatim}
\usepackage{pdfsync}
 \newtheorem{thm}{Theorem}[section]
 \newtheorem{lemma}[thm]{Lemma}
 \newtheorem{cor}[thm]{Corollary}
 \newtheorem{prop}[thm]{Proposition}
  
\theoremstyle{definition}
 
 \newtheorem{rmk}[thm]{Remark}
 \newtheorem{defn}[thm]{Definition}
 \newtheorem{claim}[thm]{Claim}

\newcommand{\N}{\mathbb{N}}

\def\L{{\mathbb{L}}}
\def\E{{\text{E}}}

\def\R{{\text{R}}}

\newcommand{\g}{\Gamma}
\def\S{{\text{S}}}
\def\F{{\text{F}}}
\def\D{{\text{D}}}
\def\L{{\text{L}}}

\def\Hom{{\text{Hom}}}
\def\U{{\text{U}}}
\def\B{{\text{B}}}
\def\C{{\textbf{C}}}
\def\P{{\mathbb{P}}}
\def\cC{{\mathcal{C}}}
\def\cB{{\mathcal{B}}}
\def\cA{{\mathcal{A}}}

\def\cD{{\mathcal{D}}}
\def\cS{{\mathcal{S}}}
\def\cG{{\mathcal{G}}}
\def\cH{{\mathcal{H}}}

\def\cI{{\mathcal{I}}}
\def\cIg{{\mathcal{I}_\text{\tiny G}}}
\def\cU{{\mathcal{U}}}
\def\cT{{\mathcal{T}}}

\def\cF{{\mathcal{F}}}
\def\cW{{\mathcal{W}}}
\def\cA{{\mathcal{A}}}
\def\cE{{\mathcal{E}}}
\def\cN{{\mathcal{N}}}
\def\cF{{\mathcal{F}}}
\def\cL{{\mathcal{L}}}

\def\ld{{\text{L}_{\mathcal{D}}}}
\def\rd{{\text{R}_{\mathcal{D}}}}
\def\fd{{\text{F}_{\mathcal{D}}}}
\def\ud{{\text{U}_{\mathcal{D}}}}
\def\ug{{\text{U}_{\mathcal{G}}}}
\def\fg{{\text{F}_{\mathcal{G}}}}
\def\uh{{\text{U}_{\mathcal{H}}}}

\def\cO{{\mathcal{O}}}
\def\o{{\small{\text{o}}}}
\def\op{{\small{\text{op}}}}

\def\Ds{{\mathcal{D}_s}}
\def\i{{\text{i}}}
\def\P{{\text{P}}}
\def\cWg{{\cW_\text{\tiny G}}}

\def\b#1{{\textbf{#1}}}
\def\F{{\text{F}}}
\def\gG{{\Gamma_\text{\tiny G}}}

\def\Hom{{\text{Hom}}}

\def\G{{\text{G}}}
\def\Ob{{\mbox{Ob}}}
\def\gl{{\text{gl}}}

\def\cTg{{\mathcal{T}_\text{\tiny G}}}

\def\map#1#2#3{\text{Map}_{#1}(#2,#3)}

\def\pig#1{\Pi_{\text{\tiny G}}^{#1}}
\def\lpig#1{\Pi_{\text{\tiny G,#1}} }
\def\gGs#1{\g_{\text{\tiny G,#1}} }

\def\GL{{\text{GL}}}

\definecolor{error}{rgb}{0.8,0,0}

\begin{document}

\title{Units of Equivariant Ring Spectra}
\author{Rekha Santhanam }
\begin{abstract}

It is well known that very special $\Gamma$-spaces and grouplike $\E_\infty$
spaces both model connective spectra. Both these models have equivariant analogues in the case when the group acting is finite. 
Shimakawa  defined the category of  equivariant $\Gamma$-spaces and showed that special equivariant  $\Gamma$-spaces determine positive equivariant spectra. Costenoble and Waner \cite{MR1012523} showed that  grouplike equivariant $\E_\infty$-spaces determine connective equivariant spectra.

We show that with suitable model category structures the category of  equivariant
$\Gamma$-spaces is Quillen equivalent to the category of equivariant $\E_\infty$
spaces. We  define the units of equivariant ring spectra in terms of equivariant
$\Gamma$-spaces and show that the  units of  an equivariant ring spectrum determines a connective equivariant spectrum.

\end{abstract}
\maketitle

\section{Introduction}
There are several space level models for the category of spectra. 
 Segal \cite{MR0353298} developed the notion of very-special $\g$-spaces to model connective spectra. May \cite{MR0339152} showed that group-like $\E_\infty$-spaces model connective spectra.

May and Thomason \cite{MR508885} gave a comparison of these models and showed that they are indeed equivalent.
However, the model theoretic viewpoint was missing and the equivariant case was not considered. 
We show that the two models of equivariant infinite loop spaces, namely, equivariant $\E_\infty$-spaces and equivariant $\g$-spaces are equivalent. 

We interpret the infinite loop space of an $\E_\infty$-equivariant ring spectrum as an equivariant $\g$-space. We then describe the units of equivariant spectra in terms of equivariant $\g$-spaces.

\subsection{Background and Results}

Let $R$ be an $\E_\infty$-ring spectrum. Then $\pi_0(R)$ defines a monoid and we can consider its unit components.
Define $\text{GL}_1 R$ to be the following pullback of spaces
\begin{equation*}
\xymatrix{\text{GL}_1 R \ar[r] \ar[d] & \Omega^{\infty} R \ar[d] \\
                     (\pi_0(R))^{\times} \ar[r] & \pi_0(R) }
\end{equation*}
May, Quinn and Ray \cite{MR0494077} showed that $\text{GL}_1(R)$ is a grouplike $\E_\infty$-space  and hence determines a connective spectrum which is denoted by $\text{gl}_1 R$.

 The theory of units of ring spectra was developed to understand the obstruction theory  \cite{MR0494077} for $\E_\infty$-orientations  on cohomology theories and to classify these orientations. 
Further the classifying space of the multiplicative units of a cohomology (ring) theory parametrize its twistings  \cite{units}, as in the case of twisted K-theory\cite{MR2172633}.

A recent result of Freed, Hopkins and Teleman  \cite{MR2365650} relates twisted equivariant K-theory  of a compact lie group with the representations of the loop group of the lie group.   Atiyah and Segal \cite{MR2172633} and Freed, Hopkins and  Teleman \cite{MR2365650} give a geometric construction of twisted equivariant K-theory.  This construction does not use homotopy theoretic methods. Further equivariant orientation theory is not as well understood as the non-equivariant case. 
We expect that the twistings of equivariant K-theory will be parametrized by the units of equivariant K-theory as in the non-equivariant case. We also hope that the units of equivariant ring spectra will give a better perspective on equivariant orientation theory.

May's machine describing equivariant infinite loop spaces via equivariant grouplike $\E_\infty$-spaces can be applied directly to construct the unit equivariant spectrum associated to the unit space of an  equivariant $\E_{\infty}$ ring-spectrum. According to May [private communication], the details have been understood in principle since the early 1980s, although the theory has still not been written up. The details of how equivariant $\E_{\infty}$-spaces describe equivariant infinite loop spaces have been discussed by Costenoble and Waner in \cite{MR1012523}. 

In this article, we give a comparison theorem, between the two models of equivariant infinite loop spaces.  We use the comparison theorem to give a construction of the unit space of equivariant $\E_\infty$-ring spectrum in terms of equivariant $\g$-spaces (Defn \ref{segalgamma}, Defn \ref{shimgamma}).  

Let $G$ be a finite group. Shimakawa defined the notion of $\Gamma_G$-spaces \cite{MR1003787} and showed that special $\Gamma_G$-spaces are equivalent to positive $G$-spectra (Defn. \ref{positive} ). We develop this notion of equivariant $\Gamma$-spaces further and show that very special $\Gamma$-spaces are equivalent to equivariant infinite loop spaces  in Theorem \ref{fixonvsp}. 

We describe a model structure on  the category of equivariant $\Gamma$-spaces where the special $\gG$-spaces are the fibrant objects.  We prove that this category is Quillen equivalent to the category of equivariant $\E_\infty$-spaces with the model structure inherited from that on the underlying category of $G$-spaces in Theorem  \ref{comparison}. We expect this equivalence will respect the symmetric monoidal structures on the categories. This is discussed in Remark \ref{symmonidal}.

If $X$ is a very special $\gG$-space then $X(\b 1)$ is a equivariant infinite loop space. Given a special $\gG$-space we show that the $G$-space represented by the orbit diagram of invertible fixed point components defines an equivariant infinite loop space cf. Lemma \ref{orbit}. 

Beginning with an equivariant $\E_\infty$-ring spectrum we define the group of units of equivariant $\E_\infty$-ring spectra as a very-special  equivariant gamma space in Definition \ref{unitsdefn}. 
Our definition of $\GL_1$ matches with  the usual notion of units of commutative ring spectra when the group action is trivial. 

In Appendix \ref{discussion}, we discuss further why our definition of equivariant units is a good analog of the non-equivariant definition.  As alluded to in Appendix  \ref{discussion} in a later paper joint with Chenghao Chu we will discuss the Quillen equivalence between the category of  equivariant $\g$-spaces and the category of equivariant spectra. There we will also discuss an equivariant analog of Segal's method of obtaining $\g$-spaces from symmetric monoidal categories.

All of our constructions are valid only when the group acting is finite. 
If $G$ is not finite then Blumberg \cite{MR2286026} shows that one cannot use the model of $\Gamma_G$-spaces. The equivariant infinite loop space theory is not as well understood when the group acting is not finite.  

There has been some work in the direction of describing equivariant infinite loop spaces in the compact Lie  group case by Caruso and Waner \cite{MR777436}. However, very little is known so far. 

\begin{rmk}
We expect that the notion of orientations arising from the equivariant space  $\GL_1$ for the  Eilenberg-Maclane spectra of  Burnside Green functors should be related the notion of equivariant orientation theory described by May, Costenoble and Waner \cite{MR1856029} for equivariant bundles when the group acting is finite.  At this point we do not have any results in this direction.
\end{rmk}

\subsection*{Acknowledgements}
A large portion of this article is my Ph.D Thesis completed at the University of Illinois, Urbana-Champaign under the guidance of Charles Rezk. 
I would like to thank Charles Rezk for his valuable advice and guidance. I also want to thank  Peter May for his suggestions and feedback on the article.

\section{Notation}
\begin{itemize}
\item  Let $\cT$ denote the category of compactly generated based topological spaces, morphisms being continuous based maps.
\item Let $\cW$ denote the category of pointed CW-complexes.
\item Let $n$, $m$, $p$ and $r$ denote natural numbers.
\item We will denote the unit of adjunction of an adjoint pair by $\eta$ and the counit by $\epsilon$. \item Denote the category of sets by $\cI$ and the category of finite $G$ sets by $\cI_\text{\tiny G}$.
\item  Let $\cC$ be any topological category and $A$ be an object of $\cC$. Then denote the corepresentable functor $\cC(A,\_ ) $ from  $\cC \to \cT$ by $\cC_A$ and representable functor  $\cC(\_ , A)$ from $ \cC \to \cT$ by $ \cC^A$. 

\end{itemize}

\section{Equivariant Infinite Loop Space machines}

Let $G$ be a compact lie group. Let $\cU$ denote the complete universe of  real representations of $G$, namely, $\cU$ is a collection of $G$-representations containing  the trivial representation and countably  many copies of irreducible representations.
\begin{defn}

 A prespectrum $X$ is a collection of $G$ spaces indexed on  finite dimensional subspaces, namely,  $V, W$ of $\cU$ with $G$-maps $S^W \wedge X_V \to X_{V\oplus W}$. If  the adjoint maps are $G$-weak equivalences then $X$ is called a $\Omega$ $G$-spectrum.

\end{defn}

For the rest of this article we will assume that $G$ is a finite group.

\subsection{Equivariant $\g$-spaces }

Shimakawa  \cite{MR1003787} constructed an equivariant analogue of $\g$-spaces. We now describe equivariant $\g$-spaces. 

Let $G\cT$ denote the category with objects based $G$-spaces and morphisms  continuous $G$-maps. 
A map of $G$-spaces $X \xrightarrow{f} Y$ is a $G$-homotopy (weak) equivalence if for every $H<G$,
$$X^H \xrightarrow{f^H} Y^H $$ is a homotopy (weak) equivalence. \\

Define $\cTg$ to be the category whose  objects are the same as that of  $G\cT$ but morphisms are all maps between based $G$-spaces.
The category $\cTg$ is enriched over $G$-spaces. Given two based $G$-spaces $X$ and $Y$, for any $f: X \to Y$ and $g\in G$ we define,
$$ g.f(x):= gf(g^{-1}x).$$
Thus, 
the space of all maps
$\cTg(X,Y)$ has  a $G$-action by conjugation.

\begin{defn}\label{segalgamma}
Let $\g$ denote the skeletal category of finite pointed sets with pointed set maps as morphisms. Denote the $n+1$ element set $\{ 0, 1, \cdots, n \}$ by $\b n$ where $0$ is the marked point. 
The category $\g$ is a topological category with discrete topology on the morphism sets. Note that our $\g$ is Segal's $\Gamma^{\op}$.\\

Define a category  $\g[\G\cT]$  to be the category whose objects are continuous functors $X$  from $\g$ to $G\cT$ such that $ X(\b 0)$ is a point. Morphisms in this category are natural transformations. 
\end{defn}

\begin{defn}\label{shimgamma}
Let $\G\g$ denote the skeletal category of finite pointed $G$-sets  (where the $G$ action preserves the marked point) with $G$-pointed maps. 
Let $\gG$ be the category with the same objects as $\G\g$ but with morphisms being all  pointed set maps, The category $\gG$ is $G$-enriched. The $G$-action on $\gG(S, T)$ is by conjugation as before.

Define the category $\gG[\cTg]$ to have objects  continuous $G$-functors  $X$ from $\gG$ to $G$-spaces such that $X(\b 0) $ is a point. We refer the objects of $\gG[\cTg]$ as equivariant $\g$-space or as $\gG$-spaces. Morphisms  in $\gG[\cTg]$ are $G$-natural transformations.

Denote the category of functors $X: G\g \to \G\cT$  such that $X(\b 0 ) = \ast $ by $ G\g[G\cT]$
\end{defn} 

Let $S$ denote a finite pointed $G$-set.
Let $p_s : S \to \b 1$ for $s \in  S$ be the morphism defined
\[
p_s(t)  = \begin{cases} 1 & \mbox{ if } s=t \\
                                                0 & \mbox{ if }   s\neq t \end{cases}
 \]
Let $X$ be a $\gG$-space. 
The projections  $p_s$ induce a map $ \theta : S \wedge X(S) \to X(\b 1) $ defined by $\theta(s,x):=X(p_s) (x)$.

Since $ \gG(S,\b1)\to \cTg(X(S),X(\b 1)) $ is a $G$-map and $g.p_{s} = p_{g^{-1}s}$ it is easy to show that the map $\theta$ is a $G$-map.

\begin{defn}\label{eqspecial}
Let $X$ be a $\gG$-space.
If the adjoint map $ X(S) \to \cTg(S,X(\b 1))$ is a $G$-weak equivalence then $X$ is defined to be a \emph{special }$\gG$-space. \\ 
\end{defn}
Define  the map $\b 2 \xrightarrow{\mu} \b 1$ to be such that $ \mu(0) = 0$ and $\mu(1)=1=\mu(2)$. Let $H$ be a subgroup of $G$ and  $X$ be a special $\gG$-space.  Then up to homotopy the map
$$( X(\b 1)^H)^2 \xleftarrow{\sim} X(\b 2)^H \xrightarrow{\mu} X(\b 1)^H $$
induces a monoidal structure on $X(\b 1)^H$.
\begin{defn}
Let $X$ be a special $\gG$-space. 
If  for every $H<G$, the space  $\pi_0 X(\b 1)^H$ is a group under the  monoid structure induced by specialness condition on $X$,  then $X$ is defined to be a \emph{very-special} $\gG$-space. \\
\end{defn}

Given any $\gG$-space $X$, the $G$-functor $X : \gG \to \cTg$ has a left Kan extension from the category of $G$-CW-complexes to $\cTg$. Denote the left Kan extension again by $X :\cWg \to \cTg$, where $\cWg$ is the $G$-enriched category of based $G$-CW complexes.
 Let $V$ and $W$ be a $G$-representations.
Then, the adjoint map to the isomorphism $ S^{V} \wedge S^W  \to  S^{V\oplus W}$ induces the following map
\begin{eqnarray*}
S^V &\to &\map{}{X(S^W)}{X(S^{V\oplus W})} \\ 
\implies S^V \wedge X(S^W)  & \to &   X(S^{V \oplus W}).
\end{eqnarray*}
Thus every $\gG$-space defines a $G$-prespectrum. Shimakawa \cite{MR1003787} shows that a special $\gG$-space defines  a positive $\Omega$-$G$-spectrum. 
\begin{defn}\label{positive}
A $G$-prespectrum $X$ is an  positive $\Omega$- $G$-spectrum if for every $G$-representation $V$ such that $V^G\neq \phi$, the map $X(V) \to \Omega^V X(W)$ is a $G$-weak equivalence.
\end{defn}

The following proposition is an important observation (due to Shimakawa and May) which we will use extensively.

\begin{prop}\label{eqcat}\cite{MR1132161}
Let $\i$ be the inclusion functor from $\g$ to $\gG$ taking sets to $G$-sets with trivial $G$-action. Then
there exist an adjoint pair of functors
\begin{equation*}
\xymatrix{ \gG[\cTg] \ar@<.5ex>[r]^{\i} & \g[G\cT] \ar@<.5ex>[l]^{\P} }
\end{equation*}
which induce an equivalence of categories.
\end{prop}
\begin{proof}
Let $X$ be a functor $\g \to \G\cT$. For any finite $G$ set, define $\gGs{S}$ to be the  $G$- functor $\gG \to \cTg$ as $\gGs{S} (T) = \gG(S,T) $ for all finite $G$-sets $T$.

 Define the functor $\P X : \gG \to \cTg$ at a $G$-set $S$ as the left Kan extension
$$\P X (S)= \gGs{S} \otimes_\g X ,$$ defined to be the coequalizer
\begin{equation*}
\xymatrix{\coprod \limits_{m,n} \gG(\b n, S)  \times \g(\b m, \b n) \times X(\b m) \ar@<.5ex>[r] \ar@<-.5ex>[r]  & \coprod \limits_m \gG(\b m, S) \times  X(\b m) \ar[r] & \P X(S), }
\end{equation*}
where one of the maps is given by the functoriality of $X$ and the other is composition in $\gG$ given via inclusion of $\g(\b m, \b n) \to \gG(\b m, \b n)$ but giving the sets trivial $G$-action. 
Let $S$ be a finite $G$-set and $f : S \xrightarrow{\cong}\b n $ as sets. The $G$-action on $S$ can be described by a group morphism $ \rho: G \to \Sigma_n$. Define $X(\b n)_\rho$ to be the $G$-space  $X(\b n)$ with the $G$-action defined as follows: Given an element $x \in X(\b n)$ and $g \in G$,
 $g x = g X(\rho(g)) x = X(\rho(g) (gx) $ since $X(\rho(g))$ is a $G$-map.

 Claim: $\P X(S) \cong X(\b n)_\rho$.
 Reason: 
 $\P X(S) = \coprod \gG(\b m, S) \times X(\b m) / \sim $, where $\sim$ is defined as follows. For any $h' \in \gG(\b m ,S)$, $h \in \gG(\b n, \b m) =\g(\b n, \b m) $ and $x \in X(\b n)$ we have,
 $(h',X(h) x)  \sim (fh,x)$.
 
Fix $f: S \to \b n$ to be an isomorphism of sets. This induces a group morphism $ \rho : G \to \Sigma_n$ such that for any $s \in S$, $$ f g(s) = \rho(g) f (s).$$

Define $X(\b n)_\rho$ as before, then we have a map $\beta : \P X(S) \to X(\b n)_\rho$ 
as $\beta(h, x) = X(fh)(x)$. This is a $G$map and is invertible with inverse 
$\theta : X(\b n)_\rho \to \P X(S)$ defined as $\theta(x)= (f^{-1},x)$.

Therefore,  $X \cong \i \P  X $.

Let $Y$ be an object in $\gG[\cTg]$. 
Let $S$ be a finite $G$ set with $|S| =n $. Then the $G$-set  is  completely described by $(\b n, \rho : G \to \Sigma_n)$  up to a set isomorphism, where $\rho$ describes the $G$-action on $S$.  
For any $ Y \in \gG[\cTg]$ we have, $ \i Y(S) (\b n) = Y(\b n) $.  Then $\P\i Y(S) = Y(\b n)_\rho$.      

The map $\epsilon : Y(S) \to \P \i Y(S) = Y(\b n)_\rho$ is induced by the isomorphism from $f : S \to \b n$,  that is,  $\epsilon(x) = Y(f) (x)$. This is a $G$-map since $\epsilon(g.x) = Y( f) (g.x) = g Y(gf) (x) $ since $Y$ is a $G$-functor. 

But, $gX(gf) (x) = g X(\rho(g)) Y(f) (x) = g. Y(f)(x) $ in $Y(\b n)_\rho$. 

Claim: $\epsilon$ is an isomorphism. \\
Reason : Define the inverse map $\alpha : Y(\b n)_\rho \to Y(S) $ as $\alpha(y) = Y(f^{-1}) y$.  Then $\alpha$ is a $G$-map. 
\begin{eqnarray*} \alpha(g y) & =& Y(f^{-1}) (gy) = Y(f^{-1} ) Y(\rho(g)) (gy)  \\ 
& = &Y(f^{-1} \rho(g) )(gy) = Y(g f^{-1}) (gy) = g. Y(f^{-1}) (gy) = g Y(f)(y).\end{eqnarray*}

Thus, these functors induce an equivalence of $G$-categories.
\end{proof}

In $\g[G\cT]$, a $\g$-space $X$ is \emph {special} if for every $H < G$ and  homomorphism $\rho : H \to \Sigma_n$ the map
$$ X(\b n)_\rho \to (X(\b 1)^n)_\rho$$
is a $H$-weak equivalence. The group $H$ acts on $X(\b 1)^n_\rho$ as follows. For any $h \in H$ and $(x_1, \cdots, x_n ) \in X(\b 1)^n$ we have 
\[ h. (x_1,\cdots, x_n) = (h.x_{\rho(h)(1)}, \cdots, h. x_{\rho(h)(n)} )\]

Shimakawa shows that \cite{MR1132161}[pg 226] this is equivalent to the condition that $\E X $ is a special $\gG$-space.  We will switch back and forth between these two notions of equivariant $\g$-spaces depending on the situation. 

\subsection{Equivariant Operads and Monads}

Costenoble and Waner \cite{MR1012523} showed that a $G$-grouplike $\E_\infty$-space is $G$ homotopy equivalent to  an equivariant infinite loop space. 

\begin{defn}
A $G$-operad $\cD$ is an operad in the category of  $G$-spaces. The spaces $\cD(n)$ have an action by $G \times \Sigma_n$  and the operad action maps are $G$-maps commuting  with the symmetric group action. 
We assume that  $\cD(0)$ is a point  (which induces the base point on $\cD(n)$ for all $n \in \N$ via the operad structure maps) and $1 \in \cD(1)$ is fixed under the action of $G$.
\end{defn}

\begin{defn}
A $\cD$-space is a based $G$-space $X$ along with $G$-maps $$\cD (n) \times X^n \to  X $$ commuting with the operad structure and the $\Sigma_n$-action.
Maps of $\cD$-spaces are maps of $G$-spaces which are compatible with the $\cD$-action. Denote the category of $\cD$-spaces by $\cD[\cTg]$.
\end{defn}

Given a based $G$-space $X$ we can construct a free $\cD$-space
\begin{eqnarray*}  F(X) & := & \coprod\limits_{n=0}^{\infty} \cD(n) \times_{\Sigma_n} X^n/ \sim 
\end{eqnarray*} 
where, the relation is defined as follows.
Let $\sigma_j : \cD(0)  \times \cD(1) \times \cdots \cD(0)   \cdots \times \cD(1) \to \cD(j-1) $ where $\cD(0)$ is in the $j$th spot and let $ i_j: X^{j-1} \to X^j$ be the map which inserts a point in the $j$th spot.  Then for any $d \in \cD(j)$ and $x \in X^{j-1}$
the relation is given by $(c, i_j(x)) \sim (\sigma_j(c), x)$.

\begin{defn}
Let $\cD$ be a $G$-operad. Then $\cD$ is a $\E_\infty$ $G$-operad if  $\cD(n)$ is a universal $ (G ,\Sigma_n)$  principal bundle for every $n \in \N$.

\end{defn}

A $G$-space, $X$ is said to be an $\E_\infty$-space if it has an $\E_\infty$ $G$-operad acting on it.  Given an $\E_\infty$-space $X$, the operad induces a monoidal structure up to homotopy on  $X^H$ for all subgroups $H< G$. Define $X$ to be $G$-grouplike if $\pi_0(X^H)$ is a group for all subgroups $H$ of $G$.

\section {Category of Equivariant Operators} 
Both, the category of grouplike $\E_\infty$-spaces and the category of very-special $\g$-spaces model infinite loop spaces. 
May and Thomason \cite{MR508885}  showed that both these approaches to infinite loop space are equivalent.  They defined the notion of "category of operators" to construct a category which can be compared to the category of $\E_\infty$-spaces and  the category of $\g$-spaces. We generalize their ideas to the equivariant setting. 

With appropriate model category structure the category $\g[\cTg]$ models the category of equivariant $\g$-spaces.
Our theorem compares the category of equivariant $\E_\infty$-spaces with  the category  $\g[\cTg]$. 

We now introduce the notion of a "category of equivariant operators". 
\begin{defn}\label{pi}
Let $\Pi$ denote the subcategory of $\g$ with morphisms
\begin{eqnarray*}
\Pi(\b m,\b n)& = & \{\phi \in\g(\b m,\b n) /\phi^{-1} (i) \text{ has at most one element for all } i > 0 \}
\end{eqnarray*}
\end{defn}
\noindent Note that $\Pi(\b m , \b 1)$ has the maps $p_i$ for all $i = 1, 2, \cdots, n$.  
\begin{defn}
Let $\G\Pi$ denote the subcategory of $\G\g$ such that
\begin{eqnarray*}
\G\Pi( S, T)& = & \{\phi \in\G\g( S, T) /\phi^{-1} (t) \text{ has at most one element for all } t \in T \}
\end{eqnarray*}
\end{defn}
\begin{defn}
Let $\pig{}$ denote the subcategory of $\gG$ such that
\begin{eqnarray*}
\pig{}( S, T)& = & \{\phi \in\gG( S, T) /\phi^{-1} (t) \text{ has at most one element for all } t \in T \}
\end{eqnarray*}
\end{defn}
\begin{defn}
Define a $\pig{}$-space to be a covariant $G$-functor from $X: \pig{} \to \cTg$ such that $X(\b 0) =\ast$.
Define the representable $\pig{}$-spaces, $\lpig{T}$ as follows
$$ \lpig{T}( S)= \pig{}(T,S)$$
Define $X$ to be a \emph{special} $\pig{}$-space if the map  $\theta $ induced by the maps $p_s$,
 $$ X(S) \to \map{}{S}{X(\b 1})$$ is a $G$ weak equivalence.

\end{defn}
Given any pointed  $G$-space $Y$, we can construct a $\pig{}$-space $\R' Y(S):=\map{}{S}{ Y}$. This defines a $\pig{}$-space. A map $\alpha : S \to T$ induces a map $\R'(\alpha): \map{}{S}{Y} \to \map{}{T}{Y} $ given by 
$$ \R'(\alpha) (f) (t) = \begin{cases}
f(\alpha^{-1} (t)) &  \mbox{ if } |\alpha^{-1}(t) | =1 \\
\ast &  \mbox{ if } |\alpha^{-1}(t) | =0 
\end{cases}
$$

\begin{lemma}
Let $\L'$ and $\R'$ be a pair of functors
\begin{equation*}
\xymatrix{ \pig{}[\cTg] \ar@<.5ex>[r]^{\L'} & \cTg \ar@<.5ex>[l]^{\R'} }
\end{equation*}
defined as $\L' X= X(\b 1)$ for $X \in \pig{}[\cTg]$ and $\R' Y( S) = \cTg(S,Y)  $ for $Y \in \cTg$. Then $\L'$ and $\R'$ are $G$-functors adjoint to each other.
\end{lemma}
\begin{proof}
It is easy to see that $\L'$ and $\R'$ are $G$-functors.
We will denote the unit of adjunction of an adjoint pair by $\eta$ and the counit by $\epsilon$.
Note $\eta: \L' \R' Y = \R' Y (\b 1) = \map{}{\b 1}{Y} \xrightarrow{\sim}  Y $. 

Now, $\R'\L' X ( S) = \cTg(S,\L' X) = \cTg( S,X(\b 1)) $.  The maps $p_s$ induce a map
 $$ X(S) \to \cTg( S ,X(\b 1))=\cTg(S,\L' X)=\R'\L' X( S)$$
 as defined before.

Therefore, the functors $\L'$ and $\R'$ are adjoint to each other.
 \end{proof}

The proof of Proposition \ref{eqcat} can be modified to show that $\pig{}[\cTg]$ and $\Pi[\G\cT]$ are equivalent categories.
The adjoint pair $\L'$ and $\R'$ factor through to give an adjoint pair  $\L :\Pi[G\cT] \to G\cT$  defined as 
$\L X = X(\b 1) $ and  $\R: G \cT \to \Pi[G \cT]$ defined as $ \R Y (\b n) = Y(\b 1)^n $. We have the following commutative diagram :

\begin{equation*}
\xymatrix{ \Pi[\G \cT] \ar@<.5ex>[dr]^-{\P}  \ar@<.5ex>[r]^-{\L} & \G \cT \ar@<.5ex>[l]^-{\R} \ar@<.5ex>[d]^{\R'} \\
& \pig{}[\cTg] \ar@<.5ex>[ul]^{\i} \ar@<.5ex>[u]^-{\L'}    }
\end{equation*}
\begin{defn}\cite[Defn~1.1]{MR508885}
Define a  category of operators $\cG$ to be a  topological category whose objects are the sets $\b n$ and with functors from  $\Pi$ to $\cG$ and $\cG$ to $\g$ such that the induced functor from $\Pi$ to $\g$ is the inclusion of $\Pi$ in $\g$. We will assume that $\cG(\b m, \b 0) = \ast$ for all $\b m  \in \text{Ob}\cG$.

A map of category of operators $\cG$ and $\cH$ is the following commutative diagram of continuous functors.
\begin{equation*}
\xymatrix{  & \cG  \ar[dr] \ar[dd]^{v} &  \\   
\Pi \ar[ur] \ar[dr] & \ar[d]   & \g \\
& \cH \ar[ur] & \\}
\end{equation*}
\end{defn}

\begin{defn}

Define a \emph{category of equivariant operators} to be a category
of operators $\cG$ enriched over $G$-spaces. Morphisms are
morphisms of category of operators which are $G$-functors. An
equivalence of category of equivariant operators is an morphism of
category of operators which induces $G$-weak equivalence on the morphism spaces.

Define a $\cG$-space $X$ to be a covariant $G$-functor from $\cG$ to $\cTg$ such that $X(\b 0) =\ast$. Denote the category of $\cG$-spaces by $\cG[\cTg]$.
\end{defn}
Note that any category of operators is enriched over $G$-spaces  via trivial $G$-action  and is  therefore a category of equivariant operators.

Let $\cH$ be a category of equivariant operators. Given any $\b n \in \mbox{Ob} \cH$, we have an object in the category $\cH[\cTg]$ defined as $\cH^{\b n} (\b m ) = \cH(\b m, \b n )$ for all  $\b m \in \Ob \cH$.
A morphism $\cG \xrightarrow{v} \cH$ of category of equivariant operators induces a morphism from $\cH[\cTg] \xrightarrow{v^\ast } \cG[\cTg]$ defined as $v^\ast Y = Y\o v$.

\begin{prop}
Let $\cG$ and $\cH$ be categories of equivariant operators.
Let $\cG \xrightarrow{v}  \cH$ be a morphism of category of equivariant operators.
Then there exists a functor $\cG[\cTg] \xrightarrow{v_\ast} \cH[\cTg]$ left adjoint to $v^\ast$.
\begin{equation*}
\xymatrix{ \cG[\cTg] \ar@<.5ex>[r]^{v_{\ast}} & \cH [\cTg] \ar@<.5ex>[l]^{v^{\ast}}. }
\end{equation*}
\end{prop}
\begin{proof}
Given a functor $\cG \xrightarrow{v} \cH $ and functor $\cG \xrightarrow{X} \cTg$, there exists a left Kan extension of $X$ to $\cH$ defined as the coend, $v_\ast X (\b n) = \cH^{\b n} \otimes_{\cG} X$ which is given by the following coequalizer in  $\G\cT$.
\begin{equation*}
\xymatrix{\coprod_{\b m ,\b k} \cH( \b m ,  \b n) \times \cG(\b k , \b m) \times X(\b k) \ar@<.5ex>[r]^-{\mu }
\ar@<-.5ex>[r]_-{v_{k,m}} & \coprod_{m} \cH(\b m, \b n ) \times X(\b m) \ar[r] & v_{\ast} X( \b n)}
\end{equation*}

 The adjointness is easy to check.
\end{proof}

Let $\cG$ be a category of equivariant operators. Then $\cG$ defines a monad on $\Pi[\G\cT]$. For any
$\Pi$-space $X$,
$$ \F_\cG X ( \b n) := \cG^{\b n} \otimes_{\Pi} X  :=\coprod \limits_m \cG(\b m, \b n) \times X(\b m) /\sim $$
where, for $f \in \cG(\b k , \b n)$, $x \in X(\b m)$ and $\pi \in \Pi(\b m, \b k)$ we have
$(f, \pi x) \sim (f \pi, x)$.

Thus $\cG[\cTg]$ is the category of $\F_\cG$-algebras over $\Pi[\G\cT]$.

Any  pointed $G$-operad  $\cD$ induces a category of equivariant operators. Define $\hat{\cD}$ to be the category with objects being the finite sets $\b n$ and the morphism space defined as
\begin{eqnarray*}
 \hat{\cD}(\b m, \b n) & := & \coprod \limits_{\phi \in \g(\b m, \b n)} \prod \limits_{1\leq j \leq n} \cD(|\phi^{-1}(j)|) .\\
\end{eqnarray*}
It follows  that the category  $\hat{\cD}$ is a category of equivariant operators. 

The category of operators $\hat{\cD}$ induces a monad and denote the free algebra functor $\fd : \Pi[\G \cT] \to \Pi[\G \cT]$ defined as  $$\fd X := \hat{\cD} X (\b n) = \coprod \limits_m \hat{\cD}(\b m, \b n) \times X(\b m) /\sim $$
where the relation is as before.

Given a $\cD$-space $X$ by construction $\R X$ is a $\hat{\cD}$-space. Denote this induced functor on $\cD$-spaces by $\rd$. We have the following square of adjoint pairs.
\begin{equation*} \label{eqsquare}
\xymatrix{ \Pi[\G\cT] \ar@<-.5ex>[d]_-{\fd} \ar@<.5ex>[r]^-{\L} & \G\cT \ar@<.5ex>[d]^-{\F} \ar@<.5ex>[l]^-{\R}  \\
                      \hat{\cD}[\cTg] \ar@<-.5ex>[u]_-{\ud}  & \cD[\cTg] \ar@<.5ex>[u]^-{\U} \ar@<.5ex>[l]^-{\rd} }  
\end{equation*}

\begin{defn}A morphism  $\phi$ in $\g(\b m, \b n)$ is said to be effective if $\phi^{-1}(0) = 0$. It is said to be ordered if it is order preserving. The set of ordered effective morphisms from $\b m $ to $\b n$ is denoted by $\cE(\b m, \b n)$.

\end{defn}
\begin{lemma} \label{gpushoutlem}\cite[Lemma 5.5]{MR508885}
Let $X$ be a $\Pi$-space. Let $\cD$ be a $G$-operad.
For $n \geq 1$ let $\F_p\hat{\cD}X(\b n)$ be the image of $\coprod \limits_{m \leq p} \hat{\cD}(\b m, \b n) \times X(\b m)/ \sim $. Then $\hat{\cD}X(\b n)$ is the union of $F_p \hat{\cD}X(\b n)$ over all $p$.

\noindent Moreover, $\F_0\hat{\cD}X(\b n) = X(\b 0) = \hat{\cD}X(\b 0) $
and $\F_p\hat{\cD}X(\b n)$ can be constructed as the following pushout of $G$-spaces;
\begin{equation}\label{pushout}
\xymatrix{ \coprod \limits_{\alpha \in \cE(\b p, \b n)}  \prod\limits_{1\leq j \leq n} \cD(|\alpha^{-1}(j)|) \times  \prod \limits_{\Sigma(\alpha )} sX(\b{p-1}) \ar[r]^-{v} \ar[d]_{i}  & \F_{p-1} \hat{\cD} X(\b n) \ar[d] \\
\coprod \limits_{\alpha \in \cE(\b p, \b n)} \prod\limits_{1\leq j\leq n} \cD(|\alpha^{-1}(j)|) \times \prod \limits_{\Sigma(\alpha)} X(\b p) \ar[r] & \F_{p} \hat{\cD}X (\b n) .
}
\end{equation}
Here  $sX(\b{p-1}) = \coprod\limits_i \sigma_i X(\b{p-1})$ for  $\sigma_i$ are the ordered effective morphisms from $\b{p-1} \to \b{p}$ and $\Sigma(\alpha) = \Sigma_{\alpha^{-1}(1)} \times \cdots \times  \Sigma_{\alpha^{-1}(n)}$ .
The morphism $v$ takes $(\alpha,c;\sigma_i x)$ to $(\alpha \sigma_i,c;x)$.
Then $$\hat{\cD}X(\b n) = \text{colim } \F_p \hat{\cD}X(\b n)$$ where the colimit is computed in the category of $G$-spaces.
\end{lemma}

\begin{lemma} \cite{MR508885} \label{identity}
Let $\cD$ be a $G$-operad.
The functor $UF$ is a monad on $G$-spaces and  $\L\ud \fd  \R = \U \F  $. In fact, $ \ud \fd \R = \R \U \F$.
\end{lemma}

By Proposition \ref{adjoint}, the functor $\rd$ has a left adjoint $\L_{\cD}$ and we have the following diagram.

\begin{equation} \label{eqsquare2}
\xymatrix{ \Pi[\G\cT] \ar@<-.5ex>[d]_-{\fd} \ar@<.5ex>[r]^-{\L} & \G\cT \ar@<.5ex>[d]^-{\F} \ar@<.5ex>[l]^-{\R}  \\
                      \hat{\cD}[\cTg] \ar@<-.5ex>[u]_-{\ud} \ar@<.5ex>[r]^{\ld}  & \cD[\cTg] \ar@<.5ex>[u]^-{\U} \ar@<.5ex>[l]^-{\rd} }
\end{equation}

\section{ Model Category Structures} \label{Gmodelcat}

We now set up the model category structure for all the categories which play a role in proving the main theorem. 

The $G$-topological category $G \cT $ admits a compactly generated model category structure where
\begin{itemize}
\item a map $X \to Y$ is a weak equivalence if  $X^H \to Y^H$ is a weak equivalence for all subgroups $H$ of $G$.
\item a map $X \to Y$  is a fibration if $X^H \to Y^H$ is a Serre fibration for all subgroups $H$ of $G$.
\item cofibrations are maps with left lifting property with respect to all trivial fibrations.
\end{itemize}
The sets  $I=\{(G /H \times \S^{n-1})_+ \to (G/H \times \D^n )_+/ H<G, n\geq1 \}$ and $J=\{(G /H \times \D^n)_+ \to( G/H \times \D^n \times  \text{I})_+  / H<G, n\geq1 \}$ are the generating cofibrations and trivial cofibrations in $G \cT$.

The following result is well known. This proof is an adaptation of the non-equivariant case.

\begin{thm}
Let $\cD$ be an pointed $G$-operad. The category of $\cD$-spaces forms a model category with weak equivalences and fibrations defined on the underlying category of  $G$-spaces.
\end{thm}
\begin{proof}
Let $\cD$ be the monad corresponding to the operad $\cD$. The category $ G\cT$ forms a cofibrantly generated model category with weak homotopy equivalences and Serre fibrations. The maps $ (G/H \times \S^{n-1})_+ \to (G/H  \times \D^{n})_+$ and $ (G/H \times \D^n )_+\to (G/H \times\D^n \times \text{I})_+$ are the generating cofibrations and acyclic cofibrations respectively. By \cite[Prop~5.13]{MR1806878} we need to show that the maps  

\noindent $\cD(G/H \times \S^{n-1})_+ \to \cD(G/H \times \D^{n}) $  and 
$ \cD(G/H \times  \D^n)_+ \to \cD(G/H \times \D^{n} \times \text I)_+ $ for $n \geq 1$, satisfy the cofibration hypothesis \cite[5.3]{MR1806878} and that the monad $\cD$ preserves reflexive coequalizers.

Reflexive coequalizers of spaces preserve finite products. Also, colimits commute with coequalizers implies $\cD$ preserves reflexive coequalizers.

Thus we need to show that 
\begin{itemize}
\item [(i)] for any $\cD$-algebra $Y$
\begin{equation*}
\xymatrix{ \cD(G/H \times  \S^{n-1} )_+\ar[r] \ar[d] & Y \ar[d] \\
                     \cD(G/H \times  \D^n )_+\ar[r] & \cD(G/H \times  \D^n)_+ \coprod_{\cD ( G/H \times \S^{n-1})_+} Y }
\end{equation*} 
the pushout is a Hurewicz cofibration . 
\item[(ii)] Every relative $\cD J$-cell complex is a weak equivalence. 
\end{itemize}

Note that 
\begin{itemize}
\item  $  Y \to \cD S^0 \coprod Y$  is a Hurewicz cofibration. 
\item  $Y\coprod_{\cD (G/H \times  \S^{n-1})_+} \cD(G/H \times \D^n)_+  = \B(Y, \cD(G/H \times  \S^n )_+, \cD( G/H \times  \ast)_+)$ and the degeneracy maps are 
Hurewicz cofibrations. 
\item  Hence $ Y \to  Y \coprod _{ \cD(G/H \times  \S^{n-1})_+} \cD(G/H \times \D^n)_+ $ is a Hurewicz cofibration.
\end{itemize}
Similar ideas can be used to show that every $\cD J$-relative cell complex is a weak equivalence.

\end{proof}
\subsection{Diagram Categories}

\begin{defn}\cite{MR1806878}
Let $\cA$ be a topological category. Let $\cA[\G\cT]$ denote the category of covariant functors from $\cA \to \G\cT$.
A map of $\cA$-spaces $X\to Y$ is said to be a level equivalence and a level fibration if for every object $a \in \cA$ and subgroup $H$ of $G$, the map $X(a)^H \to Y(a)^H $ is a weak equivalence and a Serre fibration respectively. A map of $\cA$-spaces is said to be a $q$-cofibration if it has the left lifting property with respect to all level acyclic fibrations.
A map of $\cA$-spaces $X \to Y $ is said to be an h-cofibration if $X(a)\to Y(a)$ is a Hurewicz $G$-cofibration (has $G$- homotopy extension property) for all $ a \in \text{Obj} \cA$.

\end{defn}
For every $a \in \cA$ we have an adjoint pair of functors,
\begin{equation*}
\xymatrix{ \G \cT \ar@<.5ex>[r]^{\F_a} & \cA[\G\cT] \ar@<.5ex>[l]^{\E_a} }
\end{equation*}
defined as $ \E_a (X) = X(a)$ and $ \F_a(A) (b)= \cA(a,b)\wedge A$ 

\begin{thm} \cite[Theorem~6.5]{MR1806878}\cite[Theorem~III.2.4]{MR1922205}
Let $\cA$ be a $G$-topological category. The category $\cA[\G\cT]$ admits a
 level model category structure where the weak equivalences are level equivalences, fibrations are level fibrations and cofibrations are $q$-cofibrations. Then $\cA[\G\cT]$ forms a compactly generated topological model category with the level model structure. The set of maps $\F_a I$ and $\F_a J$ for all objects $a$ of $\cA$ are the generating cofibrations and generating trivial cofibrations.
\end{thm}
\begin{proof}
The category $\cA[\G\cT]$ is complete and cocomplete since the colimits and limits are evaluated in the underlying category of $G$-spaces. In order to show that the model structure on $\G\cT$ lifts to $\cA[\G\cT]$, we need to show that the sets $\F_a I $ and $\F_a J $ satisfy the cofibration hypothesis. This follows from the adjointness of $\F_a$ and the model category structure on $\G\cT$.
\end{proof}

\begin{cor}
The category $G\Pi[\G\cT] $ is a compactly generated model category with the level model category structure.
Then the sets $\cF_{S} I$ and $\cF_{S} J$ for all $S \in \Ob G\Pi{}$ are the generating cofibrations and generating acyclic cofibrations.
\end{cor}

Let $S$ be a $G$-set. Let $\lpig{S}$ be an object of  $\pig{}[\cTg]$ defined as $\lpig{S}(T)= \pig{}(S, T)$. Then by restricting to the subcategory $G\Pi$ this also defines a $\G\Pi$-space.
The projection morphisms $p_s$ where $s \in S$ induce a map $ S  \wedge\lpig{1} \to \lpig{S}  $ in  $G\Pi[\G\cT]$.
By \cite[Thm~4.1.1]{MR1944041}  the left localization of $G\Pi[\G\cT]$ with respect to the set $V:=\{S \wedge\lpig{1} \to \lpig{S}/ S \in \Ob G\Pi \}$ exists.

Define a $G\Pi{}$-space $X$ to be a $V$-local object if for every map in $Z\to W$ in $V$,
$\map{}{W}{X} \to \map{}{Z}{X}$ is a $G$-weak equivalence. Further a morphism in  of $G \Pi[\G\cT]$ spaces $X \to Y$ is defined to be a $V$-local equivalence if for every $V$-local object $Z$, the map $\map{}{Y}{Z} \to \map{}{X}{Z}$ is a $G$-weak equivalence. \\

Then in the localized model category structure on $G\Pi[\G\cT]$
\begin{itemize}
\item weak equivalences are $V$-local equivalences,
\item cofibrations are cofibrations in the level model category structure
\item fibrations are maps with right lifting property with respect to trivial cofibrations and
\item fibrant objects are the $V$-local objects.\\
\end{itemize}

The category $\pig{}[\cTg]$ is equivalent to $\Pi [\G \cT]$.  Therefore, $\P$ is a right adjoint. 
Define a functor $ \pig{S}: \pig{\op} \to \cTg$ as $\pig{S}(T) = \pig{}(T, S)$. We can forget to $G\Pi$ to get a functor from $G\Pi^\op \to G\cT$
For any $X : G \Pi \to G\cT$ define $$ \E' X (S) = \pig{S} \otimes_{G \Pi} X$$ for all $ S \in \Ob{G \Pi}$.
It is easy to check that $\E'$ is the left adjoint to the forgetful functor  $\i: \pig{}[\cTg] \to G \Pi[G\cT]$.

This induces an adjoint pair of functors between $G\Pi[\G\cT]$  and $ \Pi[\G \cT] $. 
Define $\E : \Pi[G\cT] \to G\Pi[G\cT] $ as the right adjoint $ \E = \i \P$. 

\begin{equation*}
\xymatrix{ G \Pi[\G \cT] \ar@<.5ex>[dr]^-{\E'}  \ar@<.5ex>[r]^-{ \E'\circ \i} & \Pi[ \G \cT] \ar@<.5ex>[l]^-{\i \circ \P} \ar@<.5ex>[d]^{\P} \\
& \pig{}[\cTg] \ar@<.5ex>[ul]^{\i} \ar@<.5ex>[u]^-{\i}    }
\end{equation*}

Then the model category structure on $G\Pi[\G\cT]$ induces a model category structure on $\Pi[\G\cT]$ where
\begin{itemize}
\item A map $X \to Y $ is a weak equivalence if  $\E X \to \E Y$ is a weak equivalences in $G\Pi[\G\cT]$.
\item A map $X \to Y$ is a fibration if $\E X \to \E Y$ is a fibration in $G\Pi[\G \cT]$.
\item Cofibrations are maps with the left lifting property with respect to trivial fibrations.
\end{itemize}
Denote $\Pi[\G \cT]$ with the localized model category by $\Pi[\G \cT]_{\i V}$.  The notation is appropriate since by Claim \ref{twin models} we can also consider this  as localizing the induced model structure on $\Pi[\G\cT]$ with respect to $\i V$.\\

\begin{rmk} \label{fibrant}
The space $X$ is fibrant in $\Pi[\G \cT]$ if $\E X$ is fibrant in $G\Pi[\G \cT]$. Therefore,
the map $$ G\Pi[G\cT](\Pi_{G,S},\E X) \to G \Pi[G \cT](S \wedge \Pi_{G,\b 1}, \E X)$$  is a $G$-weak equivalence.
In particular,  $$ \Pi[G\cT](\i \Pi_{G,S}, X) \to \Pi[G\cT](S \wedge\i \Pi_{G,\b 1} , \E X)$$is a $G$-weak equivalence.

If $|S|=k$ then the $G$ action on $S$ can be described by an  isomorphism $\rho : \b k \to S $. Then  an argument similar to the proof of Proposition \ref{eqcat} shows that 
 $$ X(\b k)_\rho \to (X(\b 1)^k))_\rho $$ is a $G$-weak equivalence.

The space $X$ is fibrant if and only $ X(\b k)_\rho \to (X(\b 1)^k))_\rho $ is a $G$-equivalence for every $k$ and $\rho$. In this localized model category structure the fibrant objects in  $\Pi[\G\cT]$ are therefore, exactly the special $\Pi$-spaces.
\end{rmk}

 \begin{prop}\label{serreeqcat}
Let $\cG$ be a category of equivariant operators. Define
\begin{itemize}
\item a map of $\cG$-spaces  $X \to X'$ to be a weak equivalence
(fibration) if
 $\E \U_\cG X(\b n) \to \E \U_\cG X'(\b n)$ is a
weak equivalence (fibration) of $G\Pi$-spaces and  
\item a map of $\cG$-spaces to be
a cofibration if it has the left lifting property with respect to
all trivial fibrations.
\end{itemize}
Then $\cG[\cTg]$ forms a compactly generated model category with this structure.  The set of maps $\F_\cG \i  \F_S I $ and $\F_\cG \i\F_{S} J$ for all objects $S$ of $G\Pi$, are the generating cofibrations and generating trivial cofibrations.
\end{prop}
\begin{proof}
Colimits and limits  exist in the category of $G$-functors from $\cG \to \cTg$. 
 The sets  $\F_\cG \i \F_{S} I $ and $\F_\cG \i \F_{S} J $ satisfy the cofibration hypothesis. This follows from the fact that $\U_\cG$ preserves colimits, $\i \F_{S} I$ and $\i \F_{S} J$ satisfy the cofibration hypothesis and the functor $\F_\cG$ commutes with tensoring over spaces.
Now, apply the small object argument \cite[Lemma~5.3]{MR1806878} to prove the factorization axioms. The other axioms are easy to prove from definitions and since the model structure is  inherited from the model structure on $\Pi[\G\cT]$.
\end{proof}

\begin{rmk}[Claim \ref{twin models}]
Consider the category $\cG[\cTg]$ with the model structure inherited from the level model structure on $G\Pi$-spaces. 
Then $\cG[\cTg]$  has a localized model category structure with respect to the set $\{ S\wedge \F_{\cG} \i \lpig{\b 1} \to  \F_\cG \i \lpig{S}/ S \in \Ob G\Pi\}$ and this is equivalent to the model category structure on $\cG[\cTg]$ obtained from the underlying localized model category structure on $G \Pi[\G\cT]$.
\end{rmk}
\begin{cor}
The category $\hat{\cD}[\cTg]$ therefore, has a cofibrantly generated model category structure where
\begin{itemize}
\item a map $X \to Y$ is a weak equivalence (or fibration) if $\ud X \to \ud Y$ is a weak equivalence (or fibration) in $\Pi[\G\cT]_{\i V}$ and
\item cofibrations are maps of $\hat{\cD}$-spaces with left lifting property with respect to acyclic fibrations.
\end{itemize}
\end{cor}
\subsection{Quillen Equivalences}

\begin{prop}\label{pig}
 The adjoint functors
\begin{equation*}
\xymatrix{ \Pi[\G\cT]_{\i V} \ar@<.5ex>[r]^-{\L} & \G\cT \ar@<.5ex>[l]^-{\R} }
\end{equation*}
induce a Quillen equivalence between $\Pi[\G\cT]$ and $\G\cT$.
\end{prop}
\begin{proof}
Consider $ G\Pi[\G\cT]$ with the level model structure.  The functor $\R': \G\cT\to G\Pi[\G\cT] $ takes weak equivalences and fibrations to level weak equivalences and fibrations.
The adjoint pair $\L'$ and $\R'$ form a Quillen pair between $G\Pi[\G\cT]$ and $\G\cT$. Since the model structure on $\Pi[\G\cT]$ is induced by the model structure on $G\Pi[\G \cT]$ and we have the following diagram 
\begin{equation*}
\xymatrix{ \Pi[\G \cT] \ar@<.5ex>[dr]^-{\E}  \ar@<.5ex>[r]^-{\L} & \G \cT \ar@<.5ex>[l]^-{\R} \ar@<.5ex>[d]^{\R'} \\
& G\Pi[G \cT]. \ar@<.5ex>[ul]^{\i} \ar@<.5ex>[u]^-{\L'}    }
\end{equation*}
The adjoint pair $\L$ and $\R$ is a Quillen pair between $\Pi[\G\cT]$ and $\G\cT$.
We need to show that the adjoint pair
\begin{equation*}
\xymatrix{\Pi{}[\G\cT]_{\i V} \ar@<.5ex>[r]^-{\L} & \G\cT \ar@<.5ex>[l]^-{\R} }
\end{equation*}
induces a Quillen equivalence.  

Let $Y$ be a based $G$-space. Let $X \xrightarrow{\sim}\R' Y $ be a cofibrant replacement in $G\Pi[\G\cT]$ and therefore $\i \R' Y$ is a cofibrant replacement in $\Pi[\G\cT]$.  In particular,  $\E \i X \to \E \R' Y $ is a level $G$-weak equivalence  and $\E\i X$ is  cofibrant  in $G\Pi[\G\cT]_V$.  Then $\L X = \L' \E \i X=  X(\b 1) \to \L' \E \i \R' Y =\L \R  Y $ is a $G$-weak equivalence.

Let $X$ be a cofibrant-fibrant object in $\Pi[\G\cT]$. Then $\L X$ is cofibrant and also fibrant since all objects of $\G\cT$ are fibrant. Further $\R \L X (\b n) = \map{}{\b n}{X(\b 1)}$. Then  being fibrant in $\Pi[\G\cT]_{\i V}$ implies $ \E X \to \E \R \L X $ is a level $G$-weak equivalence and hence a weak equivalence. 

Thus $\L$ and $\R$ induce a Quillen equivalence.

\end{proof}
We would like to understand what it means for $\hat{\cD}X$  to be special. 

\begin{rmk}
Let $\cD$ be a $G$-operad. Let $\hat{\cD}$ be the category of equivariant operators induced by $\cD$. Let $X$ be a $\Pi[\G\cT]$-space $X$. For any $k$ and $\rho: G \to \Sigma_k$, we have 

\begin{eqnarray*}
(\hat{\cD}X(\b k))_\rho & = &     {\hat{\cD}^{\b k}}_\rho \otimes_\Pi X \\
& := & \coprod \limits_m  (\coprod\limits_{\phi \in \g(\b m, \b k)_\rho} \prod\limits_{1 \leq j \leq k} \cD(|\phi^{-1}(j)|)_\rho )  \times  X(\b m) /\sim ,
\end{eqnarray*}
where the relation is given as follows. Let $f \in \hat{\cD}(\b m, \b k)$, $y \in X(\b n)$ and $\alpha \in \Pi(\b n, \b k) $. Then $(f, \alpha y) \sim (f \circ \alpha, y)$.

Let $H$ be a subgroup of $G$ and $\rho : G \to \Sigma_k$ be a group homomorphism. 
The $G$ action on $\hat{\cD} (X(\b k))$ is via the $G$-action on $\g(\b m, \b k)$ for all $\b m$ and diagonal action on $\cD(|\phi^{-1}(j)|) \times X (\b m) $ for all $\phi \in \g(\b m, \b k)$.
Taking fixed points, when  $ \phi \in  \g( \b m, \b k)_\rho^H $, the map $\rho$ acts by identity on $D(|\phi^{-1}(j)|)$ for all $1\leq j \leq k$.

In fact,
\begin{eqnarray*}
(( \hat{\cD}X(\b k))_\rho)^H & = & \coprod \limits_m  \coprod \limits_{\phi \in \g(\b m, \b k)_\rho^H} \prod \limits_{1 \leq j \leq k} (\cD(|\phi^{-1}(j)|))^H  \times X(\b m)^H /\sim.
\end{eqnarray*}
where the relation $\sim$ is as defined before.
\end{rmk}

Let $X$ be an object of $\Pi[\G\cT]$. Let $H$ be a subgroup of $G$. Let $\rho: G \to \Sigma_n$ be a group homomorphism. By Lemma \ref{gpushoutlem} we have a filtration for $\hat{\cD}X(\b n)$.

\noindent Then by above remark  and since fixed points preserve pushouts, note that  $\F^H_0\hat{\cD}X(\b n)_\rho = X(\b 0)^H = \hat{\cD}X(\b 0)_\rho^H $ and $\F_p^H \hat{\cD}X(\b n) $ is the pushout  of the following diagram;
\begin{equation*}\label{Gpushout}
\xymatrix{ \prod \limits_{\alpha \in \cE(\b p, \b n)^H_\rho}  \coprod \limits_{1 \leq j \leq n} \cD(|\alpha^{-1}(j)|)^H  \times _{\Sigma(\alpha)} sX(\b{p-1})^H \ar[r]^-{v} \ar[d]_{i}  & \F_{p-1}^H \hat{\cD} X(\b n)_\rho \ar[d] \\
\coprod \limits_{\alpha \in \cE(\b p, \b n)^H_\rho} \prod \limits_{1\leq j \leq n} \cD(|\alpha^{-1}(j)|)^H \times_{\Sigma(\alpha)} X(\b p)^H \ar[r] & \F_{p}^H \hat{\cD}X (\b n)_\rho
}
\end{equation*}
Here $sX(\b{p-1})^H = \cup_i \sigma_i X(\b{p-1})^H$ where $\sigma_i$ are the ordered effective morphisms from $\b{p-1} \to \b{p}$.
The morphism $v$ takes $(\alpha,c;\sigma_i x)$ to $(\alpha \sigma_i,c;x)$.

\begin{lemma}
If $X$ is a cofibrant in $\Pi[\G\cT]$ then the map $i$ is a h-cofibration. 
\end{lemma} 
\begin{proof}
Consider the $\Pi^{\op} \times \Pi$ space $\Pi'(\b m, \b n) = \Pi (\b m ,\b{n-1})$. Then for any ordered effective morphism $\b{n-1} \to \b n$ induces a level-wise h-cofibration $\Pi' \to \Pi$ and hence $\Pi' \circ X \to \Pi \circ X$  
is a cofibration in the Hurewicz-Strom model structure by Proposition \ref{action}. This implies in particular, $sX(\b{p-1})^H \to X(\b p)^H$ is a  $\Sigma_p$-h-cofibration for all $H <G$, and $\i$ is  a  h-cofibration.
\end{proof}

\begin{lemma} \label{geq}
Let $\cD$ be a $\Sigma$-free $G$-operad, that is, $\cD(n)$ is a free $\Sigma_n$-space for all $n$.
 Let $X \to X'$ be a level $G$-weak equivalence of cofibrant objects in $\Pi[\G\cT]$. Then $ \E \hat{\cD}X \to  \E \hat{\cD} X'$ is a weak equivalence in $G \Pi[\G\cT]$.
\end{lemma}
\begin{proof}
If $X \to X'$ is a level $G$-weak equivalence then
\small{
\[
\coprod \limits_{\alpha \in \cE(\b p, \b n)^H_\rho}  \prod \limits_{1\leq j \leq n} \cD(|\alpha^{-1}(j)|)^H  \times_{\Sigma(\alpha)} sX(\b{p-1})^H  \to  \coprod \limits_{\alpha \in \cE(\b p, \b n)^H_\rho}  \prod \limits_{1\leq j \leq n} \cD(|\alpha^{-1}(j)|)^H  \times_{\Sigma(\alpha)} sX'(\b{p-1})^H \] 
}

and
\small{
\begin{eqnarray*}
\coprod \limits_{\alpha \in \cE(\b p, \b n)^H_\rho} \prod\limits_{1 \leq j \leq n} \cD(|\alpha^{-1}(j)|)^H \times_{\Sigma(\alpha)} X(\b p)^H \to \coprod\limits_{\alpha \in \cE(\b p, \b n)^H_\rho} \prod \limits_{1 \leq j \leq n} \cD(|\alpha^{-1}(j)|)^H \times_{\Sigma(\alpha)} X'(\b p)^H && 
\end{eqnarray*}  }
\normalsize{
are weak equivalences for all $H<G$.}

As in the non-equivariant case,  we can show that if $X$ is cofibrant then the map $sX(\b{p-1})^H \to X(\b p)^H $ is a $\Sigma_p$-cofibration  for all $H<G$. Since $\Sigma_k$ acts freely on $\cD(k)$, we have that $i$ is a Hurewicz-cofibration. Therefore
the pushout diagram preserves weak equivalences.

By induction, $\F_p^H \hat{\cD}X(\b n)_\rho \to \F_p^H \hat{\cD}X'(\b n)_\rho$ is a weak equivalence for all subgroups $H$ of $G$. Thus inducing a $G$-weak equivalence $\hat{\cD}X(\b n)_\rho \to \hat{\cD}X'(\b n)_\rho$, that is, a $G$-weak equivalence $\E \hat{\cD} X \to \E\hat{\cD} X$.
\end{proof}

\begin{prop}
Let $\cD$ be a $\Sigma$-free $G$-operad  and $X$ be a cofibrant-fibrant  object in $\Pi[\G\cT]$ in the localized model category. Then $\ud\fd X $  is fibrant.
\end{prop}
\begin{proof}
A $\Pi$-space, $X$ is fibrant if for every $ n \in \N$ and homomorphism $\rho : G \to \Sigma_n$, the map $ X(\b n)_\rho \to X(\b 1)^n_\rho $ is a $G$-weak equivalence.
Let $X'$ be the $\Pi$-space defined as $ X'(\b n) := X(\b 1)^n $. By Lemma \ref{geq} we have that $ \E \hat{\cD} X \to  \E \hat{\cD}  X'$ is a level $G$-weak equivalence of $\Pi$-spaces. By construction this implies $\hat{\cD}X \to \hat{\cD}X'$ is a level $G$-weak equivalence of $\hat{\cD}$-spaces.

\end{proof}

Let $X$ be an object of $\Pi[G\cT]$. Then we can define a simplicial object in $\Pi[G\cT]$ using the monad structure  on $\Pi[G \cT]$ due to $\cG$. Given an object of  $\Pi[G\cT]$, let $\B_{\ast}(\cG, \cG, Y)$ denote the simplicial object in $\Pi[G\cT]$ 
with the $n$th simplex  $\ug\fg  \cdots \ug\fg Y$ where $\ug \fg$ is applied $n$-times. The simplicial structure follows from the monad structure of $\fg\ug$. Denote the geometric realization of this simplicial object in $\Pi[G\cT]$ by $\B(\cG,\cG, Y)$.
Given a morphism of category of equivariant operators $v: \cG \to \cH$ one can similarly define $\B(\cH,\cG, Y)$.

\begin{lemma}\label{eqveq}
Let $\cG$ and $\cH$ be category of equivariant operators and $\cG \xrightarrow{v} \cH$ be a morphism of category of operators. 
 If  $\Pi \to \cG$ is a cofibration (Proposition \ref{piopspaces}) and $\E  \ug \cG^{\b m} \to  \E \uh \cH^{\b m} $ is a weak equivalence in $G\Pi[\G\cT]$ then the map of $\cH$ spaces $ \B(\cG,\cG,Y) \to \B(\cH,\cG,Y)$ is a weak equivalence.
\end{lemma}
\begin{proof}
The assumption that $\E \ug \cG^{\b m} \to \E \ug \cH^{\b m}$ is a weak equivalence implies that $ \E \ug  \B_\ast(\cG,\cG,Y) \to \E \uh \B_\ast(\cH,\cG,Y)$ is a weak equivalence of  $G\Pi$-spaces. By Proposition \ref{greedy} both these simplicial $\cH$-spaces are Reedy cofibrant  in the Hurewicz-Strom model structure. By Proposition \ref{ghurewicz} we have the lemma.
\end{proof}

\begin{thm} \label{eqcategopeq}
Let $\cG$ and $\cH$ be categories of equivariant operators.
Let  $\cG \xrightarrow{v}  \cH$ be a morphism of category of equivariant operators.
Then the following adjoint pair is a Quillen pair with the model categories inherited from the underlying category of $\Pi$-spaces.
\begin{equation*}
\xymatrix{ \cG[\cTg] \ar@<.5ex>[r]^{v_{\ast}} & \cH [\cTg] \ar@<.5ex>[l]^{v^{\ast}}. }
\end{equation*}
Further if $\Pi \to \cG$ is a cofibration (Proposition \ref{piopspaces}) and $\E  \ug \cG^{\b m} \to  \E \uh \cH^{\b m} $ is a weak equivalence in $G \Pi[G\cT]$ then the above adjoint pair form a Quillen equivalence.
\end{thm}
\begin{proof}

 Consider $\cH[\cTg]$ and $\cG[\cTg]$ with the underlying level model structure on $\Pi[G\cT]$.
Then $v^{\ast}$ takes fibrations and weak equivalences in $\cH$-spaces to fibrations and weak equivalences respectively in $\cG$-spaces since they are defined on $\Pi$-spaces. Thus $v_\ast$ and $v^\ast$ form a Quillen pair.

Note a map of $\cG$-spaces or $\cH$-spaces is a weak equivalence if it is a weak equivalence of their underlying  $\Pi{}$-spaces. In the level model structure on $\cG[\cTg]$ and $\cH[\cTg]$ all objects are fibrant. 

Let $Y$ be an cofibrant  object of $\cG[\cTg]$. By Lemma \ref{eqveq} the maps
\begin{eqnarray*}
\E \ug \B(\cG,\cG,Y) & \xrightarrow{\sim} & \E \ug Y \\
\E \uh   \B(\cG,\cG,Y) &  \xrightarrow{\sim}  & \E \uh \B(\cH,\cG,Y)  \\
\E\uh \B(\cH,\cG,Y) &\xrightarrow{\sim}& \E \uh (\cH \otimes_\cG Y) =  \E \uh v_\ast Y \\
\ug v^\ast v_\ast Y & \cong  &\uh v_\ast Y \\
\end{eqnarray*}
are $G$-weak equivalences. 
By two out of three of weak equivalences we get that implies $v^\ast v_\ast Y \to Y$ is a weak equivalence in $\cG$-spaces. 

Let $X$ be a fibrant-cofibrant $\cH$-space. Then $v^\ast X$ is a fibrant $\cH$-space. Let $Y$ such that $ Y \xrightarrow{\sim} v^\ast X$ be a cofibrant replacement $\cG$-spaces. Then $\E\uh X \xrightarrow \E \ug Y$. But we know that $$ \E \uh v_\ast Y \xrightarrow{\sim} \E \ug Y \to \xrightarrow{\sim} \E \uh X .$$

Therefore, $v_\ast$ and $v^\ast$ induce Quillen equivalence. 
\end{proof}

\begin{thm}\label{equiop}
Let $\cG$ and $\cH$ be a category of equivariant operators and $\cG \xrightarrow{v} \cH$ be a morphism of category of  equivariant operators. 
Then
\begin{equation*}
\xymatrix{ \cG[\cTg] \ar@<.5ex>[r]^{v_{\ast}} & \cH [\cTg] \ar@<.5ex>[l]^{v^{\ast}}. }
\end{equation*}
form a Quillen equivalence with the localized model category structures on $\cG[\cTg]$ and $\cH[\cTg]$.
\end{thm}
\begin{proof}
The proof follows from noting that $v_{\ast} \F_\cG = \F_\cH$ and applying Theorem 3.3.20 \cite{MR1944041} to Theorem \ref{eqcategopeq}.
\end{proof}

\begin{thm} \label{equivariantoperad}
Let  $\cD$ be a $\Sigma$-free $G$-operad such that $1 \hookrightarrow \cD(1)$ is a h-cofibration.
Let $\hat{\cD}$ denote the induced category of equivariant operators.
Then the adjoint pair of functors,
\begin{equation} \label{eqadjoint}
\xymatrix{ \hat{\cD}[\cTg] \ar@<.5ex>[r]^{\ld}  & \cD[\cTg]  \ar@<.5ex>[l]^-{\rd} }
\end{equation}
are Quillen equivalences.
\end{thm}
\begin{proof}
Let $Y \to Y'$ be fibration ( acyclic fibration) of
$\cD$-spaces. Then $\U Y \to \U Y'$ is fibration (acyclic
fibration) on the underlying category of $G$-spaces. By
Proposition \ref{pig} $\R\U Y \to \R \U Y'$ is a fibration (or
acyclic fibration). This implies $\ud \rd Y \to \ud \rd Y'$ is a
fibration (or acyclic fibration) of $\Pi[\cTg]$. Therefore the
functors $\rd$ and $\ld$ form a Quillen pair.

Note that $\rd$ creates all weak equivalences in $\hat{\cD}$-spaces.

Given a cofibrant-fibrant  $\hat{\cD}$-space $X$ the map of $\pig{}$-spaces $$ \E \ud \B(\hat{\cD},\hat{\cD}, X)= \B(\hat{\cD}, \hat{\cD}, X) \to \E\ud X$$ is a level $G$-weak equivalence.

Since $\ld$ preserves weak equivalences of cofibrant objects in $\hat{\cD}[\cTg]$, $$\U \ld \B(\fd\ud,\fd \ud X) \xrightarrow{\sim}  \U \ld  X .$$

Now, $\rd$ preserves weak equivalences of $\cD$-spaces. Hence
$$\E \rd \ld \B(\fd\ud,\fd \ud X) \to \E \rd \ld  X .$$
is a weak equivalence of $G\Pi{}$-spaces.

Therefore,
$$ \E \B(\R \U \F \L, \ud \fd, \ud X) \xrightarrow{\sim} \E \rd \ld X. $$

We have the following commutative diagram
\begin{equation*}
\xymatrix{ \E  \ud X \ar[r]   & \E \ud \rd \ld  X   \\
                      \E \B(\ud \fd, \ud \fd, \ud X)  \ar[u]^{\sim} \ar[r]^{\beta} &  \E \B( \R\U\F\L, \ud\fd, \ud X). \ar[u]^{\sim} }
\end{equation*}
where the map $\beta$ is induced by the map of triples $$ \ud \fd
\to \R \L \ud \fd \R \L \to \R\U \F \L.$$

By Lemma \ref{identity} and  Lemma \ref{pig} the map $\beta$ is a
weak equivalence of simplicial objects on $\Pi$-spaces for
cofibrant-fibrant $X$. Given a cofibrant $\hat{\cD}$-space $X$,
the bar construction is Reedy cofibrant in the Hurewicz-Strom model
category structure  by Proposition \ref{greedy}. Moreover,
geometric realization preserves weak equivalences by Proposition
\ref{ghurewicz}.
Therefore, $ \ud X \to  \ud \rd \ld X $ is a weak equivalence of $\Pi$-spaces by two out of three axiom.
Lemma \ref{quillen} implies that $\rd$ and $\ld$ induce a
Quillen equivalence.

\end{proof}
\begin{rmk} \label{symmonidal}
Lydakis \cite{MR1670245} defined a symmetric monoidal structure on the category of $\g$-spaces.  One can use this idea to define a symmetric monoidal structure on the category $\Pi[\cT]$.
The functor $R : \cT \to \Pi[\cT]$  respects the symmetric monoidal structure on $\cT$ via cartesian products. The category of operators $\hat{\cC}$  defined by an operad $\cC$ on $\cT$ defines a monad on the category $\Pi[\cT]$. By \cite{BF01220868}[Thm 2.1] it is easy to check that this monad is symmetric monoidal and hence the symmetric monoidal structure lifts  to $\hat{\cC}[\cT]$ and similarly to $\cC[\cT]$.   Then in the non equivariant case, the corresponding Theorem \ref{equivariantoperad} and Theorem \ref{equiop} respect the symmetric monoidal structures. Thus the equivalence between $\E_\infty$-spaces and $\g$-spaces  is symmetric monoidal. 
We expect this to generalize to the equivariant case. We will talk about the  monoidal structure on equivariant $\g$-spaces elsewhere. 

\end{rmk}

\section{Comparison Theorem}

Let $\cN$ denote the $G$-operad, defined as $\cN(m) = \ast$ with a trivial $G$-action. Let $\cE$ be a $\E_\infty$-$G$-operad such that $1 \to \cE(1)$ is a Hurewicz $G$-cofibration.
Then by definition for every subgroup $\Lambda < G \times \Sigma_m$ such that $\Lambda$ does not contain any non-trivial subgroups of $\Sigma_m$, $\cE(m)^\Lambda \to \cN(m)^\Lambda$ is a weak equivalence.

Consider the  category of equivariant operators induced by the operads $\cE$ and $\cN$. For any  subgroup $H$ of $G$ and $\rho: H \to \Sigma_m$, $$\hat{\cE}^{\b n }(\b m)_{\rho}= \coprod_{\phi \in \g(\b m, \b n)_\rho} \prod_{1 \leq j \leq n} \cE(|\phi^{-1}( j)|)_\rho .$$

Since $\cE $ is an $\E_\infty$-$G$-operad, $\E \hat{\cE}^{\b n}$ is  $G$-weakly equivalent to $\E \hat{\cN}^{\b n}$.

\begin{thm} \label{modelcateqgamma}
The category $\g[\G \cT]$ forms a model category with the model structure induced by level model
structure of $G\g[\G\cT]$. The localized model category of $G\g[\G\cT]$ with respect to the set  $ \{\coprod_{s\in S} \g_{G,1} \to \g_{G,S} / S  \text{ is a } G \text{set} \}$  exists and induces as model category structure on $\g[\G\cT]$ where the fibrant objects are special $\g$-spaces. This is Quillen equivalent to the localized model category structure on $\hat{\cN}[\cTg]$.
\end{thm}

\begin{thm} \label{comparison}
Let $\cE$ be an equivariant $\E_\infty$-operad such that $1 \to \cE(1)$ is a h-cofibration. 
The category of $\cE$-spaces with the model category structure induced from $G$-spaces  is Quillen equivalent to the category of $\g[\G \cT]$ with the induced localized model structure.
\end{thm}
\begin{proof}
  By Theorem \ref{equivariantoperad} we get that
 $\hat{\cE}[\cTg]$ with the localized model category structure is Quillen equivalent to $\cE[\cTg]$ with the underlying model category structure of $G$-spaces.

Since $\cE$ is an $\E_\infty$-operad for any subgroup $H$ of $G$ and group homomorphism $\rho: H \to \Sigma_n$, the space $(\cE(n)_\rho)^H$ is contractible. Note that $\hat{\cN} = \g$.
This implies that
$$ (\hat{\cE}^{\b m}(\b n)_\rho )^H \to  (\g^{\b m}(\b n)_\rho)^H =(\hat{\cN}^{\b m}(\b n)_\rho)^H$$  is a weak equivalence.
Thus $\E \cU_\cE \hat{\cE}^{\b m} \to  \E \cU_\cN \hat{\cN}^{\b m}$ is a weak equivalence of $G\Pi{}$-spaces. Theorem \ref{eqcategopeq} implies that $\hat{\cE}[\cTg]$ is Quillen equivalent to $\hat{\cN}[\cTg]$ with the localized model structures.  Hence proved.
\end{proof}

 \begin{prop} \label{main}
Let  $\cE$ be an $\E_\infty$- $G$-operad satisfying the hypothesis of Theorem \ref{eqcategopeq}. Let $X$ be an equivariant $\cE$-space. 
Then $X$ is a special equivariant $\g$-space up to a cofibrant replacement in $\g[G\cT]$.
\end{prop}
Proposition \ref{main} follows from the following Lemma and Theorem \ref{comparison}.
\begin{lemma}\label{replacement}
Let $\cD$ be a $\Sigma$-free $G$-operad and $\hat{\cD} \xrightarrow{v}\cH$ be a morphism of category of equivariant operators satisfying the hypothesis of Theorem \ref{eqcategopeq}.
Let $X$ be a fibrant-cofibrant $\hat{\cD}$-space. Then the map $$ \ud \B(\hat{\cD}, \hat{\cD},X) \xrightarrow{v} \uh \B(\cH,\hat{\cD},X)$$ is a  weak equivalence in  $\Pi[\G\cT]$ and $ \B(\cH,\hat{\cD},X)$ is a fibrant $\cH$-space in the localized model category.
\end{lemma}
\begin{proof}

By hypothesis  $\E\ud \B_\ast(\hat{\cD}, \hat{\cD},X) \xrightarrow{v} \E\uh \B_\ast(\cH,\hat{\cD},X)$ is a weak equivalence of simplicial $\Pi$-$G$-spaces. By Proposition \ref{greedy}  these spaces are Reedy cofibrant. This induces a level  $G$-equivalence of $\Pi$-spaces. This implies that $\uh \B(\cH,\hat{\cD}, X)$ is a fibrant $\Pi$-space and hence a fibrant $\cH$-space in the localized model category.
\end{proof}

\begin{proof} of  Proposition \ref{main}
 We have shown that $\cE[\cTg]$ is Quillen equivalent to $\g[G\cT]$ via a Quillen equivalence with the category $\hat{\cE}[\cTg]$. Let $X$ be an $\E_\infty$-space. Then $\R_{\cE} X$ is a fibrant $\hat{\cE}$-space. If we take its cofibrant replacement $Y \to \R_\cE X$ in the level model category structure  on $\hat{\cE} [G\cT]$ then $Y$ is fibrant-cofibrant in the localized model category on $\hat{\cE}$-spaces.
Now let $v :\hat{\cE} \to \g$ denote the morphism of category of operators induced by the contractibility of $\cE(n)$'s. 
By  Proposition \ref{replacement} the space $v_\ast \B(\hat{\cE},\hat{\cE},Y)$ is a fibrant $\g$-$G$-space.
Thus, up to a cofibrant replacement, an  equivariant $\E_\infty$-space is equivalent to a special equivariant $\g$-space.
\end{proof}

\section{$\gG$-spaces and Equivariant spectra}

Shimakawa \cite{MR1003787} generalized Segal's work to the equivariant case to show that special $\gG$-spaces model positive connective $\Omega$-$G$ spectra.
We extend Shimakawa's work to show that very-special $\gG$-spaces model connective $G$-spectra.

Given a $G$-functor $X: \gG \to \cTg$ the  left Kan extension  of $X$ from the category $\cWg$ of  based $G$-CW complexes to $\cTg$ exists. Denote the Kan  extension by $X$ again and the homotopy Kan extension of $X$ to $\cWg$ by $\tilde{X}$. 
\begin{rmk}
Note both these constructions are functorial. This defines a functor from the category of equivariant $\g$-spaces to the category of equivariant spectra. We elaborate this further in Appendix \ref{discussion}.

\end{rmk}

Let $A$ be an object of $\cWg$. Define a functor $Y_A : \g^{\op}_\text{\tiny G}\to \cTg$ as $Y_A(S) =\cTg(S,A)$. 

Then the homotopy extension is given by $\tilde{X}(A)= \B(Y_A,\gG, X)$.

Here $\B(Y_A, \gG, X)$ denotes the geometric realization of the simplicial space 

\noindent $\B_\bullet(Y_A, \gG, X)$ whose $n$ simplices are $$\coprod_{T_i} Y_A(T_n) \times \gG(T_{n-1}, T_n) \times \cdots \times \gG(T_0,T_1) \times  X(T_0),$$
face maps are compositions and degeneracy maps are defined via the natural inclusion of identity map in  $\gG(T_i,T_i)$.

The left Kan extension of $X$ at $A$ is the coequalizer
\begin{equation*}
\xymatrix{\coprod\limits_{T_0, T_1} Y_A(T_1) \times \gG(T_0, T_1) \times X(T_0) \ar@<.5ex>[r]\ar@<-.5ex>[r]  & \coprod\limits_{T_0} Y_A(T_0) \times X(T_0) \ar[r] & X(A) .}
\end{equation*}

There exists a natural map from $\tilde{X}(A)  \to  X(A)$.

\begin{lemma}\label{homtokan}
Let $X$ be an object in $\gG[\cTg]$ such that $ \i X$ is a cofibrant object in  $G\g[\G\cT]$ with the localized model category structure on $G\g$-spaces. Then  the map $\tilde{X} \to X$ is a level $G$-weak equivalence.
\end{lemma}
\begin{proof}
Let $ X$ be a representable object of $\gG[\cTg]$ denoted by $\g_{G,S}$ where $S$ is an object of $\gG$.
Then $\tilde{\g}_{G,S} \to \g_{G,S} $ is a level $G$-weak equivalence in $G\g[\G\cT]$, since all the $n$-simplices  in $\B_\bullet(Y_A,\gG,\g_{G,S})$ are degenerate for $n > 2$.

The set of maps $$I =\{\g_{G,S} \times (G /H \times \S^{n-1} )_+\to \g_{G,S} \times (G/H \times \D^n)_+ / n \in \N, S \in \text{Obj}\gG , H<G \} $$ is the set of generating cofibrations of $G\g[\G\cT]$. A cofibrant object in $G\g[\G\cT]$ can be written as a transfinite composition of maps which are pushouts of maps in $I$.

The bar construction commutes with colimits. Furthermore if the colimits are computed along cofibrations then they preserve the weak equivalences. Therefore, if $X$ is cofibrant then  $\tilde{X} \to X$ is a level $G$-weak equivalence in $G\g[\G\cT]$.

\end{proof}

\begin{thm}\cite[\text{Lem 1.4}] {MR1003787} \label{segalbit}
Let $X$ be a special $\gG$-space. Then for $G$-CW complexes  $A$ and $B$ and an object $S$ of $\gG$
\begin{itemize}
\item[(a)] the map $ \tilde{X} (S \wedge A) \xrightarrow{\tilde{\rho}} \map{}{S}{\tilde{X}(A)} $ adjoint to the evaluation map $ S \wedge \tilde{X}(S\wedge A) \xrightarrow{\rho} X(A)$ is a $G$-weak equivalence.
\item[(b)] If $\tilde{X}(A)$ is  $G$-grouplike and $A \to B$ is a $G$-cofibration then $\tilde{X} (A) \to  \tilde{X}(B) \to \tilde{ X}(A/B)$ is a $G$-fibration.
\end{itemize}

\end{thm}

\begin{thm} \label{shimakawabit}
Let $X$ be a very-special  $\gG$-space such that $\i X$  is  a cofibrant object in $G\g[\G\cT]$ with the localized model structure. Let $V$ and $ W$ be  $G$-representations  such that $V^G \neq \{ 0 \}$. Then $X (S^V) \xrightarrow{\sim} \Omega^{W} X(S^{V\oplus W})$ and $ X(S^0) \xrightarrow{\sim} \Omega X(S^1)$ are level $G$-weak equivalences.
\end{thm}
\begin{proof}
Shimakawa \cite[\text{Thm B}]{MR1003787}, shows that  for any $G$-representations $V$ and $W$ such that $V^G \neq \{0 \}$, the map $\tilde{X} (S^V) \xrightarrow{\sim} \Omega^{W} \tilde{X}(S^{V\oplus W})$ and $ \tilde{X}(S^0) \xrightarrow{\sim} \Omega \tilde{X}(S^1)$ are  level $G$-weak equivalences.

By Lemma \ref{homtokan} we have that $\tilde{X}(S^V) \to X(S^V)$ is a level $G$-weak equivalence for all representations $V$ of $G$. Since $S^V$ is cofibrant in $\G\cT$, the map $\Omega^V{X(S^W)} \to \Omega^V{X(S^W)}$ is a $G$-weak equivalence.

For any $G$-representations $V,W$ such that $V^G \neq \{0\}$ the map $$X (S^V) \xrightarrow{\sim} \Omega^{W} X(S^{V\oplus W})$$ and $$ X(S^0) \xrightarrow{\sim} \Omega X(S^1)$$ are level $G$-weak equivalences.

\end{proof}

\begin{lemma} \label{verysp}
Let $X$ be a very special  $\gG$-space and cofibrant as an object of $G\g[\G\cT]$. Let $A$ be a based $G$-CW complex. Then $X(A)$ is $G$-grouplike.
\end{lemma}
\begin{proof}
The space  $X(\b 1) $ is $G$-grouplike implies that $X(\b 1)$ has a homotopy inverse under the monoid structure. Let $S$ be finite pointed $G$-set. Then $S$ is equivalent to $\vee_i (G /H_i)_+$ as $G$-sets, for some subgroups $H_i$ of $G$.

Since $X$ is special ,
\begin{eqnarray*}
X(S) = X(\vee_i((G/H_i)_+))  = \prod_i X((G/H_i)_+) & \xrightarrow{\sim} & \prod_i \map{}{(G/H_i)_+}{X(\b 1)}  \\
& = & \prod_i X(\b 1)^{H_i}. \\
\end{eqnarray*}

Since the above equivalence commutes with the monoidal structure, $X(S)$ is $G$-grouplike.

Let $\Delta : \b 2 \to \b 1$ denote the map of finite sets which map both 1 and 2 to 1.
Then $X(S \wedge \b 1)$ being $G$-grouplike is equivalent to the map
$$X(S \wedge \b 2)^H \xrightarrow{X(\Delta) \times X(p^2_1)} X(S \wedge \b 1)^H \times X(S \wedge \b 1)^H$$ being a $G$-homotopy weak equivalence. Now taking homotopy Kan extension preserves the homotopy equivalence since geometric realizations preserve finite products.
Therefore for any $G$-CW complex $A$, 
$$ \tilde{X}(A \wedge \b 2)^H \xrightarrow{ \tilde{X}(\Delta) \times \tilde{X}(p^2_1)} \tilde{X}(A \wedge \b 1 )^H \times \tilde{X}(A\wedge \b 1)^H$$
 But by Lemma \ref{homtokan} homotopy Kan extension is weakly equivalent to Kan extension. Thus, we have that $X(A)$ is grouplike for all $G$-CW complexes $A$ if $X$ is very-special.

\end{proof}
\begin{thm} \label{fixonvsp}
Let $X$ be a very-special  $\gG$-space such that $\i X$  is  a cofibrant object in $G\g[\G\cT]$ with the localized model structure. Then $\{X(S^V)\}$ is an equivariant $\Omega$-spectrum.
\end{thm}
\begin{proof}
By Lemma \ref{verysp} and Lemma \ref{segalbit}(b), given a very-special $\gG$-space $X$ which is cofibrant in $G\g[\G\cT]$, for any $G$-representation $V$ we have that $$X(S^V) \xrightarrow{\sim} \Omega^1 X(S^{V\oplus  \mathbb{R}}).$$
But, Thm \ref{shimakawabit} says that $X(S^0) \xrightarrow{\sim} \Omega^{V} X(S^{V\oplus\mathbb{R}})$ is a $G$-weak equivalence.

Then in the following diagram for any $G$-representation $V$
\begin{equation*}
\xymatrix{ X(S^0) \ar[d] \ar[r] & \Omega X(S^1) \ar[d] \\
          \Omega^V X(S^V) \ar[r] &  \Omega^{V\oplus\mathbb{R}} X(S^{V \oplus \mathbb{R}}), }
\end{equation*}
both the horizontal arrows and the right vertical arrow are $G$-weak equivalences.
This implies that $X(S^0) \xrightarrow{\sim} \Omega^V X(S^V)$ is a  $G$- weak homotopy equivalence.

Therefore, $\{X(S^V)\}$ is an equivariant $\Omega$-spectrum.
\end{proof}

\section{$G$-spaces and Orbit Categories}

\begin{defn}
Let $G$ be a finite group. Define the orbit category of $G$ denoted $\cO(G)$ to be the category with 
\begin{itemize}
\item left cosets $G/H$ for every subgroup $H$ of $G$ as the objects
\item and, $G$-set maps as the morphisms.
\end{itemize}
The morphism set can be identified as follows $$G\cT(G/H,G/K) \cong (G/K)^H.$$
\end{defn}
\begin{defn}
Let  an $\cO(G)$-space be a  functor from $\cO(G)^{\text{op}} $ to $\cT$.
Define the category of  $\cO(G)$-spaces be the category whose objects are $\cO(G)$-spaces and morphisms are natural transformations.  Denote this category by $\cO(G)[\cT]$. 
Define the representable $\cO(G)$-space  as $$\underline{G/H}(G/K) := \cO(G)(G/K,G/H).$$
The category $\cO(G)[\cT]$  is enriched over itself.
For any two functors $W$ and $Z$ define the $\map{}{W}{Z}$ as the $\cO(G)$-space defined by the functor $ \map{}{W}{Z} (G/H) = \map{}{\underline{G/H} \times W}{Z}$.
We use the same notation for the enriched category. 
\end{defn}
Given a $G$-space $W$, we can define a $\cO(G)$-space $\Phi W$ defined as
$$\Phi W(G/H) := G\cT(G/H,W)= W^H.$$

Note that since the category $\cTg$ is $G$-enriched, it is naturally $\cO(G)$ enriched.

The model category structure on $G$-spaces is as described in Section \ref{Gmodelcat}. 
The category of $\cO(G)$-spaces has a level-model category structure where
\begin{itemize}
\item a map $W \to Z$ is a weak equivalence (fibration)  if $W(G/H) \to Z(G/H)$ is a weak-equivalences (Serre fibration) of spaces.
\item and, cofibrations are maps with the left lifting property with respect to acyclic fibrations.
\end{itemize}

\begin{prop} \cite{MR97k:55016},\cite{MR690052}\label{fixedpoints}
The functor $\Phi $ has a left adjoint $\C$  and we have an adjoint pair of enriched functors.
\begin{equation*}
\xymatrix{ \cO(G)[\cT] \ar@<.5ex>[r]^{\C}  & \cTg \ar@<.5ex>[l]^{\Phi} }
\end{equation*}
and the categories $\cO(G)[\cT] $ and $G \cT$ with the model structures described above 
\begin{equation*}
\xymatrix{ \cO(G)[\cT] \ar@<.5ex>[r]^{\C}  & \G\cT \ar@<.5ex>[l]^{\Phi} }
\end{equation*}
is a Quillen equivalence. 

\end{prop}

\section{Units of Equivariant Ring Spectra}

We construct the group of units of a special  equivariant $\g$-space and show that it is a very special equivariant $\g$-space. We use the equivalence of equivariant $\g$-spaces and  equivariant $\E_\infty$-spaces to give a construction of the units of equivariant ring spectra. 

\subsection{Units of Special $\g$-spaces}
Denote the category of sets with set maps by $\cI$ and let $\cIg$ denote the category of $G$-sets with the morphisms being set maps. The category $\cIg$ is $G$-enriched. Given a set map $K \xrightarrow{f} L$, the action of $G$ is defined via conjugation as follows:
$$ g. f(k) = g^{-1} f (g k) $$ for all $k \in K$.

\begin{defn}
Define a $\g$-set $N$ to be a functor from $\g \to \cI$ such that $N(\b 0) = \ast$. Then $N$ is a special $\g$-set if the map
$$N(\b n) \xrightarrow{\prod_i p_i} N(\b 1)^n $$ is an isomorphism.
\end{defn}

Let $N$ be a special $\g$-set.  Let $i: \b 0 \to \b 1$ be the inclusion. Then $N(\b 1)$ is a commutative monoid via the monoidal structure given by
$$ N(\b 1) \times N(\b 1) \xleftarrow{\cong} N(\b 2) \xrightarrow{N(\mu)} N(\b 1).$$
If $\g$-set $N$ is special then $N(\b m \wedge \b2) \xrightarrow{N(\text{id} \wedge \mu)} N(\b m)^2$ is an isomorphism and we have a product structure on $N(\b m \wedge \b 1)$.

We can define the group of units of $N(\b 1)$ in terms of a very special $\g$-set. Let $N'$ be a $\g$-set defined as the pullback of the following diagram:
\begin{equation*}
\xymatrix{N'(\b m) \ar@{^{(}->}[r] \ar[d] & N(\b m \wedge \b 2) \ar[d]^{N(\text{it}\wedge\mu)} \\
N(\b m \wedge \b 0) \ar@{^{(}->}[r]^{N(i)}  & N(\b m \wedge \b 1) }
\end{equation*}
Since $\b m \wedge \b 0 = \b 0$ and $N(\b 0) = \ast $ the above construction is functorial, that is, describes a $\g$-set.
By construction  $N'(\b m)$ is the pair of invertible elements of $N(\b m)$ with their inverses. 
The $\g$-set $UN$ describing the group units of $N$ is  therefore, the image of  $N'$ under the projection onto first factor namely, $$UN(\b m) := N(\text{id} \wedge p_1) N' (\b m).$$

\begin{defn}
Let $G$ be a finite group. Define a $\gG$-set to be a $G$-functor $A$ from $\gG \to \cIg$ such that $A(\b 0)= \ast$. Let $\theta: S  \wedge A(S) \to A(\b 1 )$ be as defined in Definition \ref{eqspecial}.
A $\gG$-set $A$ is \emph{special} if the adjoint of $\theta$ induces a $G$-isomorphism, that is,
$$ A(S)^H \xrightarrow{\cong} H\cT(S,A(\b 1))$$
is an isomorphism for all subgroups $H$ of $G$.
Further $A$ is \emph{very-special} if $A(\b 1)^H$ is grouplike under the induced monoid structure for all $H <G$.
\end{defn}

Given a $G$-functor $A: \gG \to \cIg$, the functor $\Phi A $ describes a functor from $\gG$ to $\cO(G)$-sets. Now the category $\gG$ is $G$-enriched. By Proposition
\ref{fixedpoints}, the category is $\cO(G)[\cI]$-enriched.

\begin{defn}
Define a $\gG$-$\cO(G)$-set to be an $\cO(G)[\cI]$- functor from $\gG$ to $\cO(G)$-sets. Given a $\gG$-set $A$, we get an  $\gG$-$\cO(G)$-set $\Phi A$.
 \end{defn}
 Note that a $\gG$-$\cO(G)$-set can be rewritten as a $\cO(G)$-$\gG$-set.

Let $A$ be a special $\gG$-set. Define a $\gG$-$\cO(G)$-set $B$ to be the following pullback of sets,
\begin{equation*}
\xymatrix{ B(S)(G/H) \ar@{^{(}->}[r] \ar[d] & \Phi A(S \wedge \b 2)(G/H) \ar[d]^{\mu}  \\
               \Phi A(S \wedge \b 0 )(G/H)   \ar@{^{(}->}[r]^{\i} & \Phi A(S \wedge \b 1) (G/H)\\ }
\end{equation*}
The construction is functorial and therefore defines a $\gG$-$\cO(G)$-set.

Define the units of $A$ to be the $\g$-$\cO(G)$-set
$$\U A(S)(G/H) :=  \Phi A (p_1) (B(S) ) (G/H).$$

Given a $\gG$-$\cO(G)$-set $B$ the projections induce the map of $\cO(G)$-sets  similar to the $\gG$-space case. 

\begin{defn}
Define a $\gG$-$\cO(G)$-set  $B$ to be \emph{special} if the map induced by projections $p_s$
$$ B(S) \to \cO(G)[\cI](\Phi S,B(\b 1))$$
is an isomorphism of $\cO(G)$-sets.
This induces a monoidal structure on $B(\b 1)(G/H)$ for all objects $G/H$ of $\cO(G)$. If $B(\b 1)(G/H)$ is grouplike for all $H<G$ then $B$ is said to be \emph{very special}.
\end{defn}

\begin{lemma}
If $A$ is a special $\gG$-set then $\U A$ is a very-special $\gG$-$\cO(G)$-set.
\end{lemma}

Let $X$ be a special $\gG$-space. Then $\pi_0 X $ is a special $\gG$-set.
Define $\U X$ as the following homotopy pullback
\begin{equation*}
\xymatrix{ UX(S)(G/H) \ar@{^{(}->}[r] \ar[d] &  X(S)(G/H) \ar[d] \\
               \U (\pi_0 X (S)(G/H))   \ar@{^{(}->}[r]^{\i} & (\pi_0 X)(S) (G/H)\\ }
\end{equation*}
By construction,  for any map $S \to T $ since $U X $ includes into $X$ we have a map from 
$UX(S)(G/H) \to X(T) (G/H) $. But the $\pi_0 U X $ is a group and this map should factor through a group of units in of $\pi_0 X$. Therefore, 
 This $UX $ is an $\gG-\cO(G)$-space by construction.

\begin{lemma}\label{orbit}
Let $X$ be a special $\cO(G)$-space. Then $\C \U X$ is a very-special $\gG$-space.
\end{lemma}
\begin{proof}
This follows from the adjointness of the $\cO(G)$-spaces and $\cTg$.
\end{proof}

\begin{defn}
Let $X$ be a special $\gG$-space. Define the units of $X$ to be the very-special $\gG$-space, $\C \U X$.
\end{defn}

\subsection{Equivariant $\E_\infty$-ring spectra}

Denote the category of $G$-spectra by $\cS_G$.

\begin{thm} The category $\cS_G$ is a topological model category with
\begin{itemize}
\item weak equivalences being $G$-weak equivalences of $G$-spectra,
\item fibrations being Serre fibrations of $G$-spectra and,
\item cofibrations being the maps of spectra with a left lifting property with respect to acyclic fibrations.
\end{itemize}
Moreover, given a continuous monad $\cC:\cS_G \to \cS_G$ such that the category $\cC[\cS_G]$ of a $\cC$-algebras  has continuous coequalizers and
satisfies the Cofibration Hypothesis, $\cS_G$ creates a topological model structure
on $\cC[\cS_G]$.
\end{thm}
\begin{proof}
The proof  is similar to the proof of the non-equivariant version \cite[Thm~VII.4.4]{MR1417719}.
\end{proof}

There exist adjoint maps from $G$-spectra to $G$-spaces.
\begin{equation*}
\xymatrix{ \cS_G \ar@<.5ex>[r]^{\Omega^{\infty}} & \cTg \ar@<.5ex>[l]^{\Sigma^{\infty}} }
\end{equation*}
Here, $\Omega^{\infty} X = X_0$ for $0$ is the trivial representation and $\Sigma^{\infty} Y$ denotes the spectrification of $ \{\Sigma^V Y \}$.

\begin{prop} Let $\cL$ denote the linear isometries  $G$-operad. Then we have a adjoint pair of functors between equivariant  $\E_\infty$-ring spectra and $\E_\infty$-spaces.
\begin{equation*}
\xymatrix{ \cL[\cS_G ]\ar@<.5ex>[r]^{\Omega^{\infty}} & \cL[\cTg] \ar@<.5ex>[l]^{\Sigma^{\infty}} }
\end{equation*}

\end{prop}

\subsection{ Defining the Units of Equivariant $\E_\infty$ Ring Spectra}

Let $R$ be a $\E_\infty$-equivariant ring spectrum. Then $\Omega^{\infty}R$ is a $\E_\infty$- ring space. There is a forgetful functor to $\cL$-spaces which forgets the additive structure on $\Omega^{\infty}R$ due to the infinite loop space structure.
By Proposition \ref{main},
as an equivariant $\cL$-space (forgetting the additive structure) $\Omega^\infty R $ is equivalent to an equivariant special $\g$-space. We  know how to construct units of a equivariant special $\g$-space. Therefore, we can make the following definition.

\begin{defn}\label{unitsdefn}
Let $R$ be an equivariant $E_\infty$-ring spectrum and $Y$ be the special $\gG$-space equivalent to the $\cL$-space, $\Omega^{\infty} R$. Define the \emph{unit equivariant spectrum } of $R$ to be the equivariant spectrum represented by the very-special $\gG$-space $\C U Y$.
\end{defn}

\appendix
\section{Adjoint Square}
\begin{thm}\cite[Thm~3.3.10]{MR771116} 
Let $\cB$ and $\cA$ be cocomplete categories
Beck's monadicity theorem states that a functor $U : \cB \to \cA $ is monadable if and only if
\begin{itemize}
\item[(i)] $U : \cB \to \cA$ has a left adjoint.
\item[(ii)] $U$ reflects isomorphisms.
\item[(iii)] $\cB$ has coequalizers of reflexive $U$-contractible coequalizer pairs and $U$ preserves them.
\end{itemize}
\end{thm}

\begin{prop} \label{adjoint}
Let $\cD$ and $\cF$  be cocomplete categories with an adjoint pair of functors
$\xymatrix{ \cD\ar@<.5ex>[r]^{\L} & \cF \ar@<.5ex>[l]^{\R}} $ such that $LR=id$. Let $ \tilde{\cD}$ and $\tilde{\cF}$ be categories with monadable functors $U_d: \tilde{\cD} \to \cD$ and $U_f:\tilde{\cF} \to \cF$. Let $F_d$ and $F_f$ be the left adjoints to $U_d$ and $U_f$ respectively with  $L U_d F_d R = U_f F_f$. Further, let there exist $\tilde{R}: \tilde{\cF} \to \tilde{\cD} $ with the  following commuting diagram of adjoint functors, namely $R U_f = U_d \tilde{R}$.
  \begin{equation}
\xymatrix{ \cD \ar@<-.5ex>[d]_{F_d} \ar@<.5ex>[r]^{L} & \cF \ar@<.5ex>[d]^{F_f} \ar@<.5ex>[l]^{R}  \\
                      \tilde{\cD} \ar@<-.5ex>[u]_{U_d}  &
\tilde{\cF} \ar@<.5ex>[u]^{U_f} \ar@<.5ex>[l]^{\tilde{R}}. }
\end{equation}
Then $\tilde{R}$ has a left adjoint such that the following diagram of adjoints commutes
  \begin{equation}
\xymatrix{ \cD \ar@<-.5ex>[d]_{F_d} \ar@<.5ex>[r]^{L} & \cF \ar@<.5ex>[d]^{F_f} \ar@<.5ex>[l]^{R}  \\
                      \tilde{\cD} \ar@<-.5ex>[u]_{U_d} \ar@<.5ex>[r]^{\tilde{L}}  &
\tilde{\cF} \ar@<.5ex>[u]^{U_f} \ar@<.5ex>[l]^{\tilde{R}}. }
\end{equation}
\end{prop}
\begin{proof}

For any $Y$ in $\tilde{\cD}$ we have a morphism $ F_d U_d Y \to Y$ due to adjointness. Also for any $Y$ in $\tilde{\cD}$ note we have
\begin{equation*}
\xymatrix{U_dF_d Y \ar[r]^-{\eta U_d F_d \eta} &   RLU_dF_dRL = RU_fF_fL }
\end{equation*}
The above map denoted by $\beta$ is in fact  a map of triples.
Further
\begin{equation*}
\xymatrix{F_fLU_dF_dY \ar[r] & F_fLRU_fF_fY=F_fU_fF_fLY \ar[r]^-{\epsilon} & F_fL Y}
 \end{equation*}
This gives an action $\alpha: F_f L U_d F_d \to F_f L$.

Define a functor $\tilde{L}: \tilde{\cD} \to \tilde{\cF}$ as follows. For any $X$ in $\cD$,
\begin{equation*}
\xymatrix{ F_f L U_d F_d U_d X \ar@<.5ex>[r]^-{ F_f L U_d \epsilon } \ar@<-.5ex>[r]_-{\alpha U_d }  & F_f L U_d X \ar[r] & \tilde{L}X }
\end{equation*}
This is a $U_f$-contractible coequalizer. By Beck's  monadicity theorem  the above coequalizer exists.

Claim: The functor $\tilde{L} $ is adjoint to $\tilde{R}$.

Reason:  Let $X$ be an object of  $\tilde{\cD}$. Then there exist maps
\begin{equation*}
\xymatrix{U_d X \ar[r]^-{\epsilon} &  R L U_d X \ar[r]^-{R\epsilon LU_dX} & R U_f  F_f L U_d X = U_d \tilde{R} F_f  L U_d X \ar[r] & U_d \tilde{R} \tilde{L} X .}
\end{equation*}
The last map is the coequalizing  map in the definition of  $\tilde{L}$. Since $U_dF_d \to RF_f U_fL$ is a map of triples we get a map of algebras
\begin{equation*}
\xymatrix{ U_d F_d U_d X  \ar[d] \ar[r] & RF_fU_f L R F_f U_fL U_d X \ar[d] \\
                             U_d X \ar[r] & R F_f U_f L U_d X }
\end{equation*}
This  gives a map from $ X \to \tilde{R} \tilde{L} X$ since $\tilde{\cD}$ is a equivalent to the category of $U_dF_d$ algebras on $\cD$.

Since $U_d$ is monadable, we can think of $\tilde{\cD}$ as space of $U_dF_d$ algebras. Then above diagram says that we have a map of $U_dF_d$ algebras $X \to \tilde{R}\tilde{L} X$.

Given any $Y$ in $ \cF$, $\tilde{L} \tilde{R} Y$  will  be given by the coequalizer
\begin{equation*}
\xymatrix{ F_f L U_d F_d U_d \tilde{R}Y \ar@<.5ex>[r]^-{ F_f L U_d \eta } \ar@<-.5ex>[r]_-{\alpha U_d }  & F_f L U_d \tilde{R} Y \ar[r] & \tilde{L} \tilde{R}Y}
\end{equation*}

Using the fact that $U_d \tilde{R} = R U_f$ and $LR=id$ we get that the coequalizer diagram is

\begin{equation*}
\xymatrix{ F_f L U_d F_d R U_f Y \ar@<.5ex>[r]^-{ F_f L U_d \eta } \ar@<-.5ex>[r]_-{\alpha U_d }  & F_f U_f Y \ar[r] & \tilde{L} \tilde{R}Y}
\end{equation*}

This reduces to
\begin{equation*}
\xymatrix{ F_f U_f F_f  U_f Y \ar@<.5ex>[r]^-{ F_f U_f \eta } \ar@<-.5ex>[r]_-{\eta F_f U_f }  & F_f U_f Y \ar[r] & \tilde{L} \tilde{R}Y}
\end{equation*}

But this gives an isomorphism  $\tilde{L} \tilde{R} Y \to Y$.

Thus we have an adjoint pair.
\end{proof}

Let $\cA$ and $\cB$ be model categories.
A functor $\U: \cB \to \cA$ creates weak equivalences in $\cA$ if a map $B \to B'$ is a weak equivalence in $\cB$ if and only if $\U B \to \U B'$ is a weak equivalence.
\begin{lemma} \label{quillen} \cite[lemma~A.2]{MR1806878}
Let $\U: \cB \to \cA$ and $\F: \cA \to \cB$ be a Quillen adjoint pair. Then $(\U,\F)$ form a Quillen equivalence if
$\U$ creates weak equivalences in $\cB$ and for all cofibrant objects $A$ of $\cA$, the map $A \to \U\F A$ is a weak equivalence in $\cA$.
\end{lemma}

\begin{lemma} \label{twin models}
Let $\cD$ be a model category and $S$ be a set of maps in $\cD$ such that the localization of $\cD$ with respect to $S$ exists and is denoted by $\Ds$. Let $T$ be a monad on $\cD$ and $\cD_T$ denote the category of $T$-algebras. Then $\cD_T$ has a model category structure inherited by both $\cD$ and $\Ds$. Localizing $\cD_T$ with respect to $TS$ we get another model category structure on $\cD_T$ and let us denote this model category by $\cD_{TS}$ for notational convenience. 

Then the model categories $\Ds_T$ and $\cD_{TS}$ are Quillen equivalent.
\end{lemma}
\begin{proof}:  

\noindent In the localized model category $\Ds$
\begin{itemize}
\item fibrant objects are $S$-local objects, namely, fibrant objects $X$ in $\cD$ such that for every morphism $Y \to Y'$ in $S$, the map $\map{}{Y'}{X} \to \map{}{Y}{X}$ is a weak equivalence.
\item weak equivalences are $S$-local equivalences, namely, morphisms $Z \to W$ such that $\map{}{W } {X} \to \map{}{Z}{X}$ is a weak equivalence for all fibrant $X$.
\item cofibrations are maps which are cofibrations  in $\cD$.
\item fibrations are maps with right lifting property with respect to acyclic cofibrations.
\end{itemize}
In the model category $\Ds_T$
\begin{itemize}
\item Weak equivalences and fibrations same as those in $\Ds$. 
 \item Cofibrations are the ones with left lifting property with respect to acyclic fibrations.
\end{itemize}
In the model category $\cD_{TS}$ 
\begin{itemize}
\item Weak equivalences are $TS$ local equivalences.
\item Cofibrations are on the underlying category $\cD$.
\item Fibrations have the right lifting property with respect to acyclic cofibrations.
\end{itemize}
Note that the free functor $F_T$ on $\cD$ is left adjoint to the forgetful functor $U_T$ on $\cD_T$. 
Thus a $TS$-local object in $\cD_T$ is exactly a $T$-algebra whose underlying space is an $S$-local object of $\cD$. Both model categories have the same fibrant objects. 
For similar reasons, both model categories have the same weak equivalences. Moreover, fibrations in ${\Ds}_T$ are fibrations in $\cD_{TS}$.  Thus one can show that the identity functors will actually induce an Quillen equivalence between the two model categories.
\end{proof}

\clearpage
\section{Cofibrant Objects} \label{appb}
\begin{prop}
The category of $G$-topological spaces forms a model category with
\begin{itemize}
\item weak equivalences as $G$-homotopy equivalences of  $G$-spaces,
\item cofibrations are Hurewicz $G$-cofibrations denoted by h-cofibrations and
\item fibrations are maps with right lifting property with respect to trivial cofibrations denoted by h-fibrations. In particular,  fibrations are Hurewicz fibrations.
\end{itemize}
\end{prop}
We will call this the Hurewicz-Strom model structure on $G$-spaces.

\begin{prop}\label{ghurewicz}
Let $X$ be simplicial object in $G$-topological spaces with the Hurewicz-Strom model structure. Then geometric realization preserves weak homotopy equivalences between Reedy cofibrant objects.
\end{prop}

\begin{lemma}\cite[I.$\delta$6.5]{MR1417719}
Let $A \to B $ be a h-cofibration of $G$-topological spaces. Then cobase change along a weak homotopy equivalence is a weak homotopy equivalence.
\end{lemma}

\begin{lemma}\label{hcofweq}
Let the following be a pushout diagram of $G$-spaces ;
\begin{equation*}
\xymatrix{ A \ar[r]^{f} \ar[d]_i & C \ar[d] \\
                    B \ar[r] & B\cup_A C.}
 \end{equation*}
If $i$ is a h-cofibration then the pushout is preserved under weak homotopy equivalences.
\end{lemma}

For rest of this section we will assume that $\G\cT$ has the Hurewicz-Storm model structure.
Consider $G\Pi[G\cT]$ with the level model category structure. Then $G\Pi[G\cT]$ is a topological model category induced by the Hurewicz-Strom model structure on $G\cT$.  Then $\Pi[G\cT]$ is a topological model category with the model structure induced by the functor $\E$.
Let $\cG$ be a category of operators and $\cG[\cTg]$ have the model category from the underlying structure on $\Pi[G\cT]$. Denote the category of simplicial objects in the $\cG[\cTg]$ by $s.\cG[\cTg]$. Consider the category $s.\cG[\cTg]$ with the Reedy model structure induced by the Hurewicz-Strom model structure on $\cG[\cTg]$.
\begin{prop}\label{greedy}
Let $\cG$ and $\cH$ be category of operators with a morphism $v: \cG \to \cH$. Let $\Pi \to \cG$ be a cofibration of $\Pi^{\op} \times \Pi$-spaces.
Let $X$ be a cofibrant $\cG$-space. Then the bar construction $\B_\bullet(\cH,\cG,X)$ is Reedy cofibrant as a simplicial object in $\Pi[G\cT]$.
\end{prop}
 Note in the case that this Proposition is applied we assume that $X$ is cofibrant in model category
 described in Section \ref{serreeqcat}.

In order to prove Proposition \ref{greedy} we reformulate the proof of a similar result by Rezk[Thesis].

Consider the category of covariant functors from $\Pi^\op \times \Pi$ to $G \cT$ denoted by $(\Pi^\op \times \Pi)[\G \cT]$.
We can define a monoidal structure on $\Pi^\op \times \Pi$-spaces as follows. For any  $\Pi^{\op} \times \Pi$-spaces $\cA$ and $\cB$ define a $\Pi^\op \times \Pi$-space as 
$$ \cA \o \cB (\b n , \b m) := \cA^{\b m} \otimes_\Pi \cB_{\b n} $$
Note $\Pi \o \cA = \cA$ and $\cA \o \Pi = \cA$.

The category of $\Pi^{\op} \times \Pi$-spaces acts on the category of $\Pi$-spaces. Let $X$ be a $\Pi$-space and $\cA$ be a $\Pi^{\op} \times \Pi$-space. Then define a $\Pi$-space $\cA(X)$ as follows
$$ \cA(X) (\b n) := \cA^{\b n} \otimes_\Pi X .$$

This defines right closed action of $\Pi^{\op} \times \Pi$-spaces on $\Pi$-spaces. Let $X$ and $Y$ be $\Pi$-spaces.
Define $$\Hom(X,Y) (\b m, \b n) :=  \cTg(X(\b m),Y(\b n) ).$$  Note this is a $\Pi^{\op} \times \Pi$ space.

Let $X$ and $Y$ be $\Pi$-spaces and $\cA$ be a $\Pi^{\op} \times \Pi$-space. Then
$$\Hom(\cA X, Y ) \cong \Hom(\cA, \Hom(X,Y)).$$

Further given  $\Pi^{\op} \times \Pi$-spaces $\cA$ and $\cB$ we get a function $\Pi^{\op} \times \Pi$-space defined is the coequalizer
\begin{equation*}
\xymatrix{ \F(\cA,\cB) (\b m, \b n)  & \ar[l]   \prod_{k} \cTg(\cA(\b k, \b m)\cB(\b k, \b n))
& \ar@<-.5ex>[l] \ar@<.5ex>[l] \prod_{k \to k'} \cTg(\cA(\b k,\b m),\cB(\b k', \b n))  \\
}
\end{equation*}

\begin{prop}
The category of $\Pi^{\op} \times \Pi$-spaces has a right closed monoidal structure and
for $\Pi ^{\op} \times \Pi$-spaces $\cA,\cB$ and $\cG$
$$ \Hom(\cA \o \cB, \cG) \cong \Hom(\cA, \F(\cB,\cG) ).$$
\end{prop}

\begin{prop}\label{piopspaces}
Define a morphism $ \cA \to \cB$ in  $(\Pi^{\op}\times \Pi)[G\cT]$
\begin{itemize}
\item to be a weak equivalence (or fibration) if $\cA(\b n,\hspace{0.1cm}  ) \to \cB(\b n,\hspace{0.1cm}  ) $ is a weak equivalence (or fibration) in $\Pi[G\cT]$ and,
\item to be a cofibration if it has the left lifting property with respect to all trivial fibrations. 
\end{itemize}

\end{prop}

\begin{prop}\label{action}
The action  of $\Pi^{\op} \times \Pi$-spaces on $\Pi$-spaces is compatible with the level model category structure of $\Pi$-spaces.
\end{prop}
\begin{proof}
Let $i: X \to Y$ be a cofibration of $\Pi$-spaces and $p:Z \to W$ be a fibration of $\Pi^\op \times \Pi$-spaces. Then we need to show that the induced maps
\begin{eqnarray*}
f:  \Hom(\ast, Z) & \to & \Hom(\ast, W)  \text{ and} \\
g: \Hom(Y,Z)  & \to &  \Hom(X, Z)  \times_{\Hom(X,W)}  \Hom(Y, W)  
\end{eqnarray*}
are fibrations.  Further we need to show that if $i$ is also a weak equivalence, then $g$ is a trivial fibration. If $p$ is also a weak equivalence, then $f$ and $g$ are trivial fibrations. 

We can reduce this to a similar diagram in $\Pi[G \cT]$ using adjointness of $\E$ and $\i$. 
The result follows from the fact that $\Pi[G\cT]$ is a topological model category. 

\end{proof}

\begin{prop}\label{product}
The monoidal structure of $\Pi^{\op} \times \Pi$-spaces is compatible with the model category structure of $\Pi^{\op} \times \Pi$-spaces.
\end{prop}
\begin{proof}
We need to show that if $i:\cA \to \cB$ is a cofibration  and $p:\cG \to \cH$ is a fibration of $\Pi^{\op} \times \Pi$-spaces, then the induced maps
\begin{eqnarray*}
\F( \ast, \cG) & \to &  \F( \ast, \cH)  \text{ and} \\
\Hom(\cB, \cG ) &  \to & \Hom(\cA, \cG) \times_{\Hom(\cA,\cH)} \Hom (\cB, \cH) 
 \end{eqnarray*}
are fibrations in $(\Pi^{\op} \times \Pi)[G \cT]$. If $i$ is also a weak equivalence then the second map is a trivial fibration. If both $i$ and $p$ are weak equivalences then both the maps are trivial  fibrations.

In order to prove the above result we need to know that if $\cA \xrightarrow{i} \cB$ is a cofibration then $\cA^{\b n} \to \cB^{\b n}$ is a cofibration of $\Pi$-spaces. This follows from the fact that fibrations of $\Pi^{\op}\times \Pi[G \cT]$ are defined level-wise.

The rest of the proof is similar to that of the previous Proposition.

\end{proof}

In order to prove Proposition \ref{greedy} we follow the proof of Proposition 3.7.3 in Rezk[Thesis].
In order to  show that $\B_\bullet(\cH,\cG,X)$ is Reedy cofibrant for a cofibrant $\cG$-space $X$,  we need to show that
$L_{n-1}\B_\bullet(\cH,\cG,X) \to \B_n(\cH,\cG,X)$ is a cofibration of $\cG$-spaces.

Let $\Pi \xrightarrow{i} \cG$ be the natural map. Then we have maps $\cG^{\otimes m} \xrightarrow{s_j} \cG^{\otimes {m+1}}$ given by $s_j= \text{id} \otimes \cdots i \otimes \cdots \otimes \text{id}$ where $i$ is in the $j$th spot and $s_0= i \otimes \text{id} \otimes \cdots \otimes \text{id}$.
Now, define $\cA_m$ to be the following coequalizer.
\[\xymatrix{\coprod_{0 \leq r <j< m-1}  \cG^{\circ {m-1}} \ar@<.5ex>[r]^{s_r} \ar@<-.5ex>[r]_{s_j}  & \coprod_{0\leq k \leq m} \cG^{\circ{m}} \ar[r] & \cA_m }
\]

There exist maps $s_k$ from $\cG^{\circ{m}} \to \cG^{\circ{m+1}} $ giving rise to a map $a: \cA_m \to \cG^{\circ{m+1}}$.

\begin{lemma}
The following diagram is a pushout square in $\Pi^{\op} \times \Pi$-spaces.
\begin{equation*}
\xymatrix{ \cA_m \circ \Pi  \ar[r]^{\text{id}\circ i}\ar[d]^{a\circ \text{id}} & \cA_m \circ \cG \ar[d] \\
                      \cG^{\circ {m+1}} \circ \Pi \ar[r] & \cA_{m+1} }
\end{equation*}
\end{lemma}
\begin{proof}
The functor $_ \circ \cG$ preserves colimits in $\Pi^{op} \times \Pi$-spaces as they are computed in the underlying category of spaces. The proof follows similar to the proof of Lemma 3.7.8 Rezk[Thesis].  
\end{proof}

\begin{lemma}
Let $\cG$ be a $\Pi^\op \times \Pi$ space such that $\Pi \to \cG$ is a cofibration of $\Pi^{\op}\times \Pi$-spaces. Then the map $\cA_{m+1} \to \cG^{\circ{n+1}} $ is a cofibration $\Pi^{\op}\times \Pi$-spaces.
\end{lemma}
\begin{proof}
Proof is by induction. By hypothesis $(\cA_0= \Pi) \to \cG $ is a cofibration. Let $\cA_{m-1} \to \cG^{\circ m}$ be a cofibration. Then since $\Pi \to \cG$ is a cofibration, by previous lemma and Proposition \ref{action} we get that $\cA_{m} \to \cG^{\circ m+1}$ is a cofibration.
\end{proof}
\begin{rmk}
Note that if $\cD$ is a well pointed operad, that is, $\ast \to \cD(1)$ is a h-cofibration. Then we can show that $\Pi \to \hat{\cD}$ is a cofibration of $\Pi^{\op}\times \Pi$-spaces. 
\end{rmk}

\begin{lemma}
Let $\cG \to \cH$ be a map of $\Pi^\op \times \Pi $-spaces, $\Pi \to \cG$ be a  cofibration of $\Pi^{\op} \times \Pi$-spaces and $X$ be a cofibrant $\Pi$-space. Then $L_n\B_\bullet(\cH,\cG,X) \to \B_n (\cH,\cG,X)$ is a cofibration in $\Pi[G\cT]$.
\end{lemma}
\begin{proof}
By previous lemma and  Proposition \ref{action}  the map $\cA_n \circ X \to \cG^{\circ n+1 } \circ X$ is a cofibration. Now, $L_n \B_\bullet(\cH,\cG,X) \cong \cH \circ \cA^{\circ n} \circ X $.  This implies from Proposition \ref{product} that $L_n \B_\bullet(\cH,\cG,X) \to \B_n(\cH,\cG,X)$ is a level cofibration.
\end{proof}

\begin{proof} of \textbf{Proposition \ref{greedy}} follows from the previous lemma.
\end{proof}

 \section{Discussion on units of equivariant ring spectra}\label{discussion}
  

Following is a discussion regarding the equivariant $\gl_1$ functor from equivariant ring spectra to equivariant spectra. Let $\cS_G$ denote the $G$ enriched category of $G$-spectra and $G \cS$ denote the category of $G$-spectra without the enrichment.
There exists Quillen pair of functors between equivariant ring spectra and equivariant $E_\infty$-spaces 
\begin{eqnarray} 
\xymatrix{ \cE[\cTg] \ar@<-.5ex>[r]_{\Sigma^\infty} & \cE[ \cS_G]\ar@<-.5ex>[l]_{\Omega^\infty} }
\end{eqnarray}
induced by the adjoint pair between $G$-spaces and $G$-spectra.  This induces an adjoint pair between the homotopy categories of $\cE[\cTg]$ and $\cE[\cS_G]$. By the results in this article, since the homotopy categories of $\cE[\cTg]$  and $\g[G\cT]$ are equivalent we have an adjoint pair between the homotopy categories of  $\g[G \cT]$ and $\cE[\cS_G]$. 

There are two relevant model structures on the category of equivariant $\g$-spaces. The one described 
in this paper is such that the fibrant objects in the category are special  equivariant $\g$-spaces, which we will denote by $\g[G \cT]_s$. 
There is a different model structure in which fibrant objects are very special equivariant $\g$-spaces which we denote by $\g[G\cT]_{vs}$. In a later paper (joint with Chenghao Chu), we show that there is a Quillen pair between the category of  equivariant $\g$-spaces and a suitable category of equivariant spectra that induces a equivalence between the homotopy category of connective equivariant spectra and homotopy category of equivariant $\g$-spaces.  We will have a Quillen pair as follows, where $\cA$ and $\cB$ denote the equivariant analogs of functors $\cA$ and $\cB$ defined by Segal \cite{MR0353298}[Prop 3.3]
\begin{eqnarray} \label{equiv}
\xymatrix{ G\cS\ar@<-.5ex>[r]_{\cA} & \g[G\cT]_{vs} \ar@<-.5ex>[l]_{\cB} }
\end{eqnarray}

Consider the functor $\text{Units}$ obtained by taking fibrant replacement in $\g[G\cT]_s$ and then applying $\GL_1$ construction to it. On the level of homotopy categories this induces a functor which is right adjoint to the identity functor of equivariant $\g$-spaces. More precisely we have a pair of adjoint functors, 

\begin{eqnarray*}
\xymatrix{ \text{ho.}\g[G\cT]_{vs}\ar@<.5ex>[r]^{Id }  & \text{ho.} \g[G\cT]_{s} \ar@<.5ex>[l]^{\text{Units}}}
\end{eqnarray*}

 

Assembling all these diagrams and noting that the Quillen pair \ref{equiv} induces an equivalence on the homotopy category of connective spectra. We can define the functor on the homotopy cateories  $\text{gl}_1 : \text{ho.} \g[G\cS] \to  \text{ho. connective } G \cS \subset  \text{ho.} G\cS$ adjoint to the functor $\cA \Omega^\infty:  \text{ho. connective } G \cS \subset\text{ho.} G \cS \to  \text{ho.} \g[G \cS]$.

\begin{rmk}
 We expect that the notion of equivariant $\g$-spaces can be extended to the notion of equivariant $\g$-spectra and one can generalize the result in this paper to a Quillen equivalence between the category of equivariant $\E_\infty$-spectra and equivariant $\g$-spectra.  Following the notation in this article, we will have, 
 
\begin{claim}
Let $\cE$ denote a $E_\infty $-$G$-operad. Then with appropriate model structures where the fibrant objects in $\g[G\cS]$ are special objects, we get a zigzag of Quillen equivalences between 
$\cE[\cS_G]$ and $\g[G\S]$
\end{claim}

We can reiterate the definition of $\gl_1$ in the equivariant case using the above claim. 
\end{rmk}

\bibliographystyle{amsplain}
\bibliography{refer}
\end{document}